\documentclass[11pt]{article}
\usepackage[utf8]{inputenc}
\usepackage{amsmath,amssymb,amscd,epsfig,bbm}
\usepackage{stmaryrd,mathabx}
\usepackage{comment}
\usepackage{color}
\usepackage[T1]{fontenc}

\usepackage{interval}
\usepackage{combelow}

\usepackage[textsize=scriptsize,textwidth=2cm]{todonotes}
\usepackage{enumitem}
\usepackage{varwidth}
\setlist{nolistsep}
\usepackage[colorlinks]{hyperref}

\usepackage{mathtools}
\makeatletter
\DeclareRobustCommand\widecheck[1]{{\mathpalette\@widecheck{#1}}}
\def\@widecheck#1#2{%
    \setbox\z@\hbox{\m@th$#1#2$}%
    \setbox\tw@\hbox{\m@th$#1%
       \widehat{%
          \vrule\@width\z@\@height\ht\z@
          \vrule\@height\z@\@width\wd\z@}$}%
    \dp\tw@-\ht\z@
    \@tempdima\ht\z@ \advance\@tempdima2\ht\tw@ \divide\@tempdima\thr@@
    \setbox\tw@\hbox{%
       \raise\@tempdima\hbox{\scalebox{1}[-1]{\lower\@tempdima\box
\tw@}}}%
    {\ooalign{\box\tw@ \cr \box\z@}}}
\makeatother

\usepackage{amsthm}

\newtheorem{defi}{Definition}[section]
\newtheorem{prop}[defi]{Proposition}
\newtheorem{theo}[defi]{Theorem}
\newtheorem{theofr}[defi]{Théorème}
\newtheorem{conj}[defi]{Conjecture}
\newtheorem{lemm}[defi]{Lemma}
\newtheorem{lemmfr}[defi]{Lemme}
\newtheorem{coro}[defi]{Corollary}

\theoremstyle{definition}
\newtheorem{rema}[defi]{Remark}
\newtheorem{exem}[defi]{Example}
\newtheorem{exems}[defi]{Examples}

\newcommand{\bdefi}{\begin{defi}}
\newcommand{\edefi}{\end{defi}}
\newcommand{\bprop}{\begin{prop}}
\newcommand{\eprop}{\end{prop}}
\newcommand{\btheo}{\begin{theo}}
\newcommand{\etheo}{\end{theo}}
\newcommand{\btheofr}{\begin{theofr}}
\newcommand{\etheofr}{\end{theofr}}
\newcommand{\blemm}{\begin{lemm}}
\newcommand{\elemm}{\end{lemm}}
\newcommand{\blemmfr}{\begin{lemmfr}}
\newcommand{\elemmfr}{\end{lemmfr}}
\newcommand{\brema}{\begin{rema}}
\newcommand{\erema}{\end{rema}}
\newcommand{\bexer}{\begin{exem}}
\newcommand{\eexer}{\end{exem}}
\newcommand{\bexems}{\begin{exems}}
\newcommand{\eexems}{\end{exems}}
\newcommand{\bconj}{\begin{conj}}
\newcommand{\econj}{\end{conj}}
\newcommand{\bcoro}{\begin{coro}}
\newcommand{\ecoro}{\end{coro}}
\newcommand{\dem}{\noindent{\bf Proof. }}

\usepackage{mathrsfs}
\renewcommand\mathcal{\mathscr}

\newcommand{\A}{{\cal A}}

\newcommand{\D}{{\cal D}}
\newcommand{\E}{{\cal E}}
\newcommand{\F}{{\cal F}}
\newcommand{\G}{{\cal G}}

\newcommand{\N}{{\cal N}}
\newcommand{\OOO}{{\cal O}}

\newcommand{\maths}[1]{{\mathbb #1}}  

\newcommand{\CC}{\maths{C}}

\newcommand{\FF}{\maths{F}}
\newcommand{\HH}{\maths{H}}

\newcommand{\NN}{\maths{N}}

\newcommand{\PP}{\maths{P}}
\newcommand{\QQ}{\maths{Q}}
\newcommand{\RR}{\maths{R}}
\newcommand{\SSS}{\maths{S}}

\newcommand{\XX}{\maths{X}}

\newcommand{\ZZ}{\maths{Z}}

\newcommand{\mmm}{{\mathfrak m}}

\newcommand{\ra}{\rightarrow}
\newcommand{\bs}{\backslash}

\newcommand{\ov}[1]{{\overline #1}} 
\newcommand{\wt}[1]{{\widetilde{#1}}}
\newcommand{\wh}[1]{{\widehat{#1}}}

\newcommand{\ga}{\gamma}
\newcommand{\Ga}{\Gamma}

\newcommand{\cqfd}{\hfill$\Box$}

\newcommand{\arcosh}{\operatorname{argcosh}}

\newcommand{\Ax}{\operatorname{Ax}}

\newcommand{\bigO}{\operatorname{O}}

\newcommand{\card}{{\operatorname{Card}}}
\newcommand{\CAT}{\operatorname{CAT}}

\newcommand{\gengeod}{\operatorname{\widecheck{\G\,}\!\!}}

\newcommand{\id}{\operatorname{id}}
\renewcommand{\Im}{{\operatorname{Im}}}
\newcommand{\Isom}{\operatorname{Isom}}
\newcommand{\Leb}{\operatorname{Leb}}

\newcommand{\mult}{\operatorname{mult}}

\newcommand\Perp{\operatorname{Perp}}

\renewcommand{\Re}{{\operatorname{Re}}}
\newcommand{\restrict}{\!|_}

\newcommand{\ssm}{\!\smallsetminus\!}
\newcommand{\stab}{\operatorname{Stab}}

\newcommand{\Vol}{\operatorname{Vol}}
\newcommand{\vol}{\operatorname{vol}}

\newcommand{\hdr}{{\HH}^2_\RR}
\newcommand{\htr}{{\HH}^3_\RR}

\newcommand{\hnr}{{\HH}^n_\RR}
\newcommand{\hnc}{{\HH}^n_\CC}

\newcommand{\hdc}{{\HH}^2_\CC}

\newcommand{\PSL}{\operatorname{PSL}}

\newcommand{\GL}{\operatorname{GL}}

\newcommand{\PGL}{\operatorname{PGL}}

\newcommand{\Sim}[2][1]{
  \mathrel{\overset{{\scriptstyle #2}}{\scalebox{#1}[1]{$\sim$}}}
}
\newcommand{\SimDiv}[1]{\overset{{\scriptscriptstyle #1}}{\displaystyle\sim}}

\newcommand{\flow}[1]{{{\tt g}^{#1}}}  

\newcommand\normalout{\partial^1_{+}}

\newcommand\normalpm{\partial^1_{\pm}}

\newcounter{fig}

\title{Counting and equidistribution of strongly reversible  \\ closed 
geodesics
  in negative curvature} \author{ Jouni Parkkonen \and
  Frédéric Paulin} \date{\today}

\begin{document}
\bibliographystyle{../alphanum}
\maketitle

\begin{abstract} 
Let $M$ be a pinched negatively curved Riemannian orbifold, whose
fundamental group has torsion of order $2$. Generalizing results of
Sarnak and Erlandsson-Souto for constant curvature oriented surfaces,
and with very different techniques, we give an asymptotic counting
result on the number of strongly reversible periodic orbits of the
geodesic flow in $M$, and prove their equidistribution towards the
Bowen-Margulis measure. The result is proved in the more general
setting with weights coming from thermodynamic formalism, and also in
the analogous setting of graphs of groups with $2$-torsion.  We give
new examples in real hyperbolic Coxeter groups, complex hyperbolic
orbifolds and graphs of groups.
\end{abstract}

\section{Introduction}
\label{sec:intro}

Let $\wt M$ be a complete simply connected Riemannian manifold with
pinched negative sectional curvature at most $-1$. Let $\Ga$ be a
nonelementary discrete subgroup of $\Isom(X)$, having {\em
  involutions}, that is, elements of order $2$. A loxodromic element
of $\Ga$ is {\em strongly reversible} if it is conjugated to its
inverse by an involution of $\Ga$. A {\it strongly reversible closed
  geodesic} in the Riemannian orbifold $M=\Ga\bs \wt M$ is the image
of the translation axis of a strongly reversible loxodromic element of
$\Ga$.  See \cite{OFaSho15} for a general discussion and an extensive
review of reversibility phenomena in dynamical systems and group
theory.

In this paper, we give an asymptotic counting and equidistribution
result of strongly reversible closed geodesics (with multiplicities
and weights) of length at most $T\ra+\infty$, generalising results of
\cite{Sarnak07,ErlSou24}.  We refer to Section \ref{sec:dyn} for the
definition of the weights, that come from the thermodynamic formalism
of equilibrium states, see for instance \cite{Ruelle04, PauPolSha15}.
In this Introduction, we restrict to the case where all weights are
equal to $1$.
 
Strongly reversible closed geodesics appear for example in
\cite{Sarnak07}, where strongly reversible loxodromic elements of
$\ga\in\PSL_2(\ZZ)$ are called {\em reciprocal elements} because of
their connections with the reciprocal integral binary quadratic forms
of Gauss.  See also \cite{BoPaPoZa14,BouKon19, BasSuz22,BasSuz25} for
recent work on reciprocal elements of $\PSL_2(\ZZ)$.  The same
terminology is used for strongly reversible elements of Hecke triangle
groups in \cite{DasGon24}, and for those in any lattice of
$\PSL_2(\RR)$ that contains involutions in \cite{ErlSou24}.  See
Corollary \ref{coro:sarnak} and Example \ref{ex:sarnak}, where we
relate our results with \cite[Thm.~2 (13)]{Sarnak07} and
\cite[Thm.~1.1]{ErlSou24}.

In order to state a simplified version of our counting and
equidistribution result, we introduce the measures that come into
play, refering to Section \ref{sec:dyn} for precise definitions, and
to \cite{BroParPau19} for more explanations and for historical
references. We denote by $\|\mu\|$ the total mass of a measure $\mu$.
We refer to Section \ref{sec:multi} for the definition of the
multiplicity of a strongly reversible closed geodesic. For instance,
the multiplicity of a primitive strongly reversible closed geodesic is
$2$, when $\wt M$ has dimension $2$ and the involutions in $\Ga$ only
have isolated fixed points in $\wt M$.

Let $\delta_\Ga$ be the critical exponent of $\Ga$. Let
$(\mu_{x})_{x\in \wt M}$ be a Patterson density for $\Ga$ and let
$m_{\rm BM}$ be the associated Bowen-Margulis measure on $T^1M=\Ga\bs
T^1\wt M$.  When $m_{\rm BM}$ is finite, then $\frac{m_{\rm BM}}
{\|m_{\rm BM}\|}$ is the unique measure of maximal entropy for the
geodesic flow on $T^1M$, see \cite{OtaPei04,DilTho25}. When $\wt M$ is
a symmetric space and $\Ga$ has finite covolume, then $\mu_{x}$ is (up
to a scalar multiple) the unique probability measure on
$\partial_\infty \wt M$ invariant under the stabiliser of $x$ in the
isometry group of $\wt M$, and $m_{\rm BM}$ is the Liouville measure,
which is then finite and mixing.  Given a nonempty proper totally
geodesic submanifold $D$ of $\wt M$, we denote by $\nu^1 D$ its unit
normal bundle, and by $\wt\sigma^+_{D}$ (resp.~$\wt\sigma^-_{D}$) the
outer (resp.~inner) skinning measure on $\nu^1 D$ for $\Ga$, which is
the pull-back of the Patterson density by the map sending a normal
vector to $D$ to the point at $+\infty$ (resp.~$-\infty$) of the
geodesic line it defines.

Let $I_\Ga$ be the set of involutions of $\Ga$ that we assume to be
nonempty.  For every $\alpha\in I_\Ga$, let $F_\alpha$ be its fixed
point set in $\wt M$.  The group $\Ga$ acts on $I_\Ga$ by conjugation.
Let $I$ and $J$ be fixed $\Ga$-invariant nonempty subsets of
$I_\Ga$. A loxodromic element $\ga$ of $\Ga$, and its associated
(oriented) closed geodesic in $M$, is {\it $\{I,J\}$-reversible} if
$\ga$ is conjugated to its inverse by an element $\alpha$ of $I$ such
that $\ga\alpha=\alpha\ga^{-1}\in J$. Let $\N_{I,J}(T)$ be the number
of $\{I,J\}$-reversible closed geodesics (counted with multiplicities,
as defined Section \ref{sec:multi}) of length at most $T$.

The $\Ga$-invariant measure $\sum_{\alpha\in I}
\wt\sigma^\pm_{F_\alpha}$ on $T^1\wt M$ induces a locally finite
measure $\sigma^\pm_{I}$ on $T^1M$, called the {\em skinning measure}
of the family $\F_I=(F_\alpha)_{\alpha\in I}$ in $T^1 M$.  When $\wt
M$ is a symmetric space and $\Ga$ has finite covolume, the measure
$\sigma^\pm_{I}$ is nonzero. It is finite except when the fixed point
set $F_\alpha$ of some $\alpha\in I$ is a geodesic line with
noncompact image in $M$. We refer to \cite{ParPau25a, ParPauSay25c}
and the end of Example \ref{ex:extendedmodular} for further
information on the surprising counting phenomena and growth that occur
when the skinning measures are infinite.

\btheo\label{theo:mainintro} Assume that the Bowen-Margulis measure
$m_{\rm BM}$ is finite and mixing for the geodesic flow on $ T^1M$,
and that the skinning measures $\sigma^+_{I}$ and $\sigma^-_{J}$ are
finite nonzero.

\hypertarget{mainintro1}{(1)} As $T\ra+\infty$, we have
\[
\N_{I,J}(T)\sim \frac{\|\sigma^+_{I}\|\;\|\sigma^-_{J}\|}
  {\delta_\Ga\,\|m_{\rm BM}\|}
\;\exp\big(\,\frac{\delta_\Ga}{2}\;T\,\big)\,.
\]

\hypertarget{mainintro2}{(2)} When $J=I$, the sum (with
multiplicities) of the Lebesgue measures along the
$\{I,I\}$-reversible closed geodesics of length at most $T$,
normalized to be a probability measure and lifted to $T^1M$, weak-star
converges on $T^1M$ to the Bowen-Margulis measure $m_{\rm BM}$,
normalized to be a probability measure.
\etheo

The proof of Theorem \ref{theo:mainintro} relates strongly reversible
closed geodesics to common perpendiculars between the fixed point sets
of involutions of $\Ga$, as explained in Section
\ref{sec:generalities}, and then uses the counting and
equidistribution results of \cite{ParPau17ETDS, BroParPau19,
  ParPau25b}, with subtle work on multiplicities in the new
averaging arguments.

\noindent\begin{minipage}{6.5cm}
\begin{center}
\includegraphics[width=3cm]{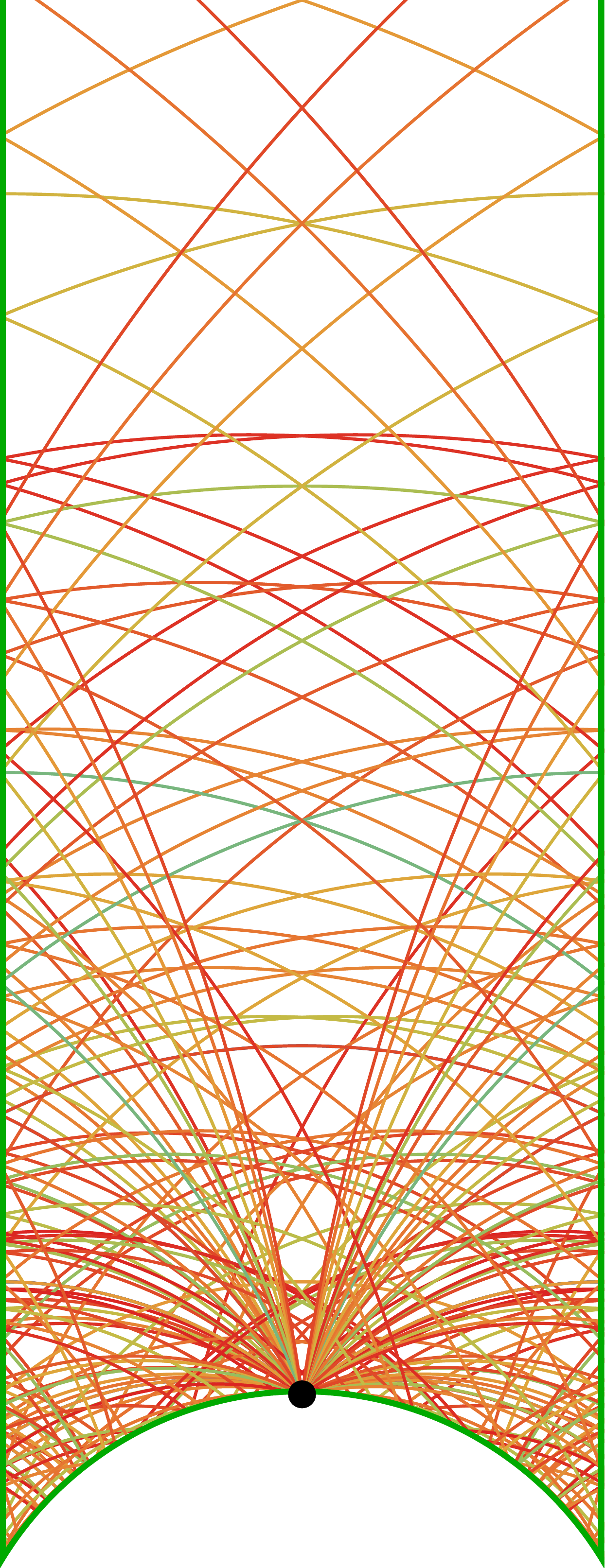}
\includegraphics[width=3cm]{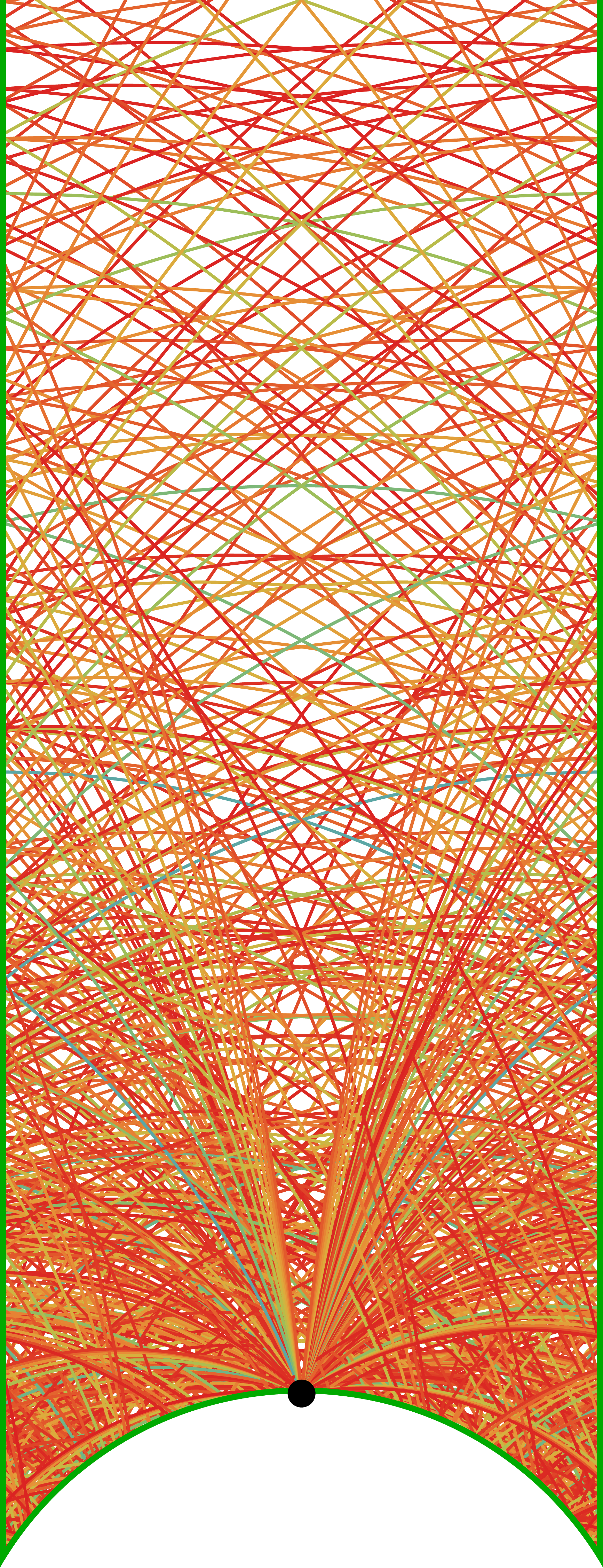}
\end{center}
\end{minipage}
\label{page:tallfigure}
\begin{minipage}{8.2cm}
\setlength\parindent{15pt} Let $\hdr$ be the upper halfspace model of
the real hyperbolic plane, let $\Ga_6$ be the Hecke triangle group of
signature $(2,6,\infty)$ and let $I$ be the set of conjugates in
$\Ga_6$ of the involution $z\mapsto-\frac 1z$.  The figure on the left (resp.~right)
shows the $\{I,I\}$-reversible closed geodesics of $\Ga_6\bs\hdr$ of
length at most $11$ (resp.~$13$) restricted to the low part of the standard
fundamental polygon of $\Ga_6$. 
See Example \ref{ex:sarnak} for more
information.

The collection of periodic orbits considered in Theorem
\ref{theo:mainintro} (\hyperlink{mainintro2}{2}) is a strict subset of
the collection of all periodic orbits (known to equidistribute to the
Bowen-Margulis measure by results of \cite{Bowen72} and
\cite{Roblin03}). The collection of $\{I,I\}$-reversible closed
geodesics of length at most $T$ grows at a rate $c\,
e^{\frac{\delta_\Ga}{2}T}$, which is considerably smaller than the
growth $c'\,T^{-1}\, e^{\delta_\Ga\,T}$ of the set of closed
geodesics of length at most $T$, where $c,c'$ are constants.
\end{minipage}

\medskip
The equidistribution of reciprocal geodesics on $\PSL_2(\ZZ)\bs\hdr$
was conjectured by Sarnak in \cite{Sarnak07}, and proved for lattices
of $\PSL_2(\RR)$ that contain involutions by Erlandsson and Souto in
\cite{ErlSou24}.  When $\wt M$ is a symmetric space and $\Ga$ is an
arithmetic lattice, we furthermore have error terms in both the
counting and equidistribution statements of Theorem
\ref{theo:mainintro} (see Sections \ref{sec:proofs} and
\ref{sec:eqproof}). The constant $\frac{\|\sigma^+_{I}\|\,
  \|\sigma^-_{J}\|}{2\,\delta_\Ga\, \|m_{\rm BM}\|}$ may be made
explicit in these cases (see Section \ref{ex:plicit}). We recover
Theorems 1.1 and 1.4 of \cite{ErlSou24} for the particular case of the
real hyperbolic plane $\wt M=\HH^2_\RR$ and $I=J=I_\Ga$ in a synthetic
way, adding an error term to their result.

We conclude this introduction with applications of Theorem
\ref{theo:mainintro}, refering to Section \ref{ex:plicit} for more
examples, in particular to Subsection \ref{subsec:complexhyp} for a
study of strongly reversible closed geodesics in complex hyperbolic
reflexion groups, as examplified by Deraux's example \ref{ex:Deraux}.

\bcoro\label{coro:Coxeter} Let $(W,S)$ be a real hyperbolic Coxeter
reflexion system in dimension $n\geq 2$ with a finite volume Coxeter
polyhedron $P$ that is compact if $n=2$.  Let $I_S$ be the set of the
conjugates by elements of $W$ of the elements of $S$.  Then there
exists $\kappa>0$ such that, as $T\ra+\infty$, we have
\[
\frac{1}{2}\N_{I_S,I_S}(T)=\frac{\Vol(\partial P)^2}
     {(n-1)\,2^n\Vol(\SSS^{n-1})\Vol(P)} \;e^{\frac{n-1}2\,T}
\big(1+\bigO(e^{-\kappa T})\big)\,.
\]
\ecoro

Corollary \ref{coro:Coxeter} is not applicable if $n=2$ and $P$ is not
compact, see \cite{ParPau25a} and Example \ref{ex:extendedmodular} for
further information.  For example,
\[
P=\big\{(z,t)\in\htr=\CC\times\RR_+: 
0\le\Re\;z \leq \frac{1}{2}\,,\ 0\le \Im\;z \leq \frac{1}{2}\
\textrm{and}\ |z|^2+t^2\geq 1\big\}
\]
is a finite volume Coxeter polyhedron in the upper halfspace model
$\htr$ of the real hyperbolic $3$-space.  The dihedral angles between
two vertical faces of $P$ are $\pi/2$ and the dihedral angles of the
edges of the spherical face of $P$ are $\pi/3$. The following figure
shows some of the boundaries at infinity of the walls of the Coxeter
system $(W,S)$ generated by reflexions in the (codimension $1$) faces
of $P$. See Example \ref{ex:Gauss} for further information.
\begin{center}
\includegraphics[width=13cm]{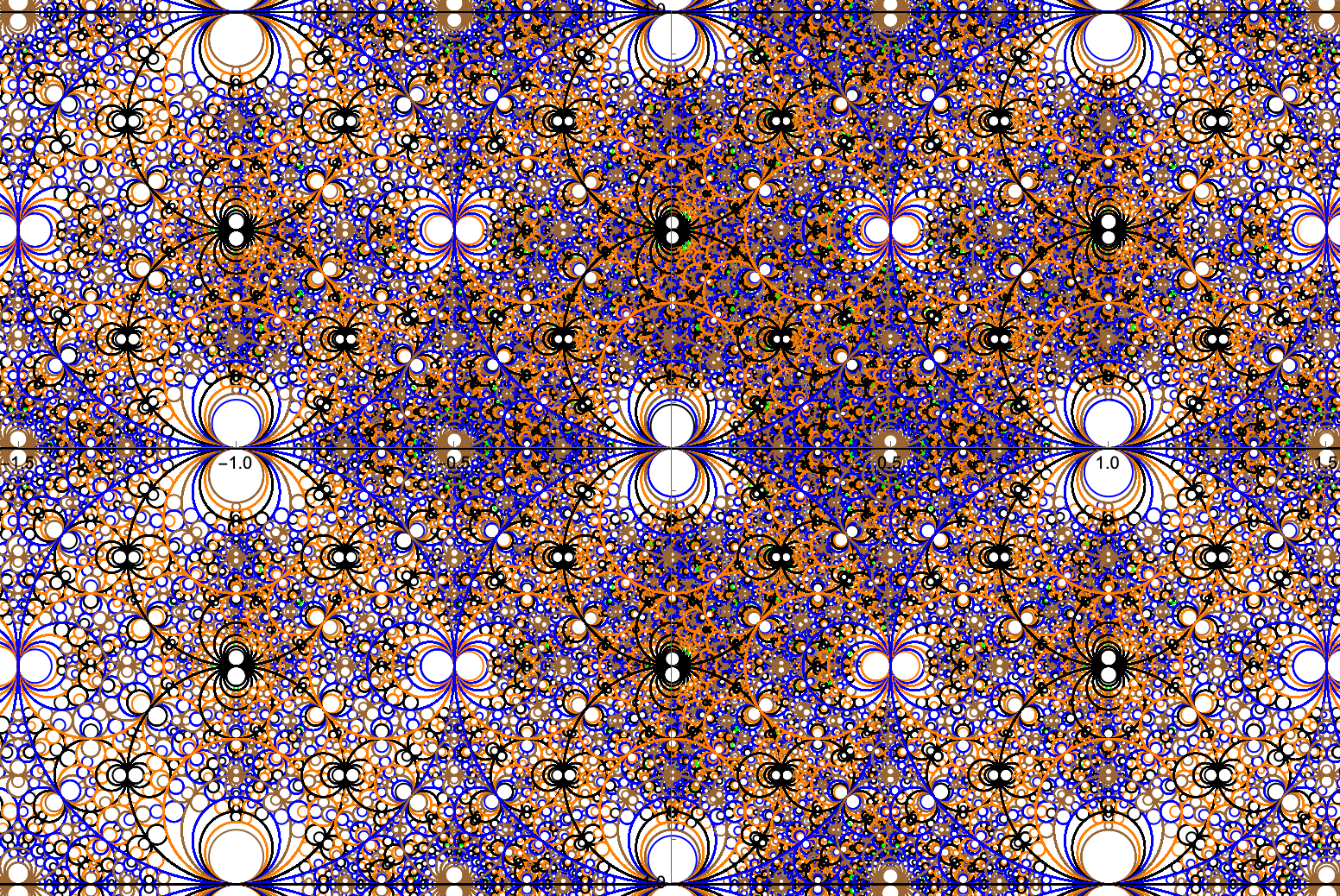}
\end{center}

In Sections \ref{sec:proofs} and \ref{sec:eqproof}, we prove more
general versions of the counting and equidistribution results stated
in Theorem \ref{theo:mainintro}, with potentials coming from the
thermodynamic formalism of equilibrium states, see \cite{PauPolSha15},
including a version of Theorem \ref{theo:mainintro} for simplicial
trees.  An application of the case of trees is given by the following
result. We refer to Subsection \ref{subsec:tree} for the proof of this
result, and for the relevant definitions.

\bcoro\label{cor:introNagao} Let $q$ be a prime power with $q\equiv
3\bmod 4$. Let $\Ga=\PGL_2(\FF_q[Y])$ and let $I_\alpha$ be the
conjugacy class in $\Ga$ of the involution $\alpha =\begin{bsmallmatrix}
0&1\\-1&0\end{bsmallmatrix}$.  The number $\N_{I_\alpha,I_\alpha}(n)$ of
con\-ju\-gacy classes of $\{I_\alpha, I_\alpha\}$-reversible
loxodromic elements of $\Ga$ whose translation length on the
Bruhat-Tits tree of $(\PGL_2, \FF_q((Y^{-1})))$ is at most $n$
satisfies, as $n\in 4\ZZ$ tends to $+\infty$,
\[
\N_{I_\alpha,I_\alpha}(n)\sim \frac{q^{\frac n2-1}}{q^2-1}\,.
\]
\ecoro

\medskip
\noindent{\small {\it Acknowledgements: } We thank the organizers of
  the conference ``Groups, Geometry and Dynamics'', June 6-10, 2022,
  in the beautiful by-the-sea CNRS conference center in Cargese, where
  the birth of this paper took place. In particular, we thank Viveka
  Erlandsson for the inspiration to write this paper, and for all the
  discussions in Cargese and later on: this paper would not exist
  without her. We thank Juan Souto for discussions. We thank Martin
  Deraux a lot for the discussions concerning Subsection
  \ref{subsec:complexhyp} and for writing \cite{Deraux24} that
  provided most of the computations needed for our Example
  \ref{ex:Deraux}.  We thank the French-Finnish CNRS IEA PaCap and the
  Magnus Ehrnooth foundation for their support. }

\section{Strongly reversible elements and common perpendiculars}
\label{sec:generalities}

Let $X$ be either a complete simply connected Riemannian manifold $\wt
M$ with dimension at least $2$ and with pinched negative sectional
curvature $-a^2\leq K\leq -1$, or the geometric realisation of a
simplicial tree $\XX$ without terminal vertices and with uniformly
bounded degrees of vertices.  We denote by $X\cup\partial_\infty X$
the geometric compactification of $X$, where $\partial_\infty X$ is
the boundary at infinity of $X$.  For every closed convex subset $A$
of $X$, we denote by $\overline{A}$ its closure in
$X\cup\partial_\infty X$.

When $X$ is a manifold, we denote by $\Isom(X)$ its locally compact
full isometry group.  When $X$ is a tree, we denote by $\Isom(X)$ the
locally compact group of automorphisms of $\XX$ without edge
inversion.  Let $\Ga$ be a nonelementary discrete subgroup of
$\Isom(X)$.  An element in $\Isom(X)$ of order $2$ is called an {\em
  involution}.  We assume from now on that $\Ga$ contains involutions.
See for instance \cite{BriHae99} for background on $\CAT(-1)$ spaces
and their discrete groups of isometries, and \cite{Serre83} for
background on group actions on trees, as well as
\cite[Chap.~2]{BroParPau19}.

For every $\ga\in\Isom(X)$, we denote by
\[
\lambda(\ga)=\inf_{x\in X} d(x,\ga x)
\]
the {\em translation length} of $\ga$. The element $\ga$ is {\it
  loxodromic} if $\lambda(\ga)>0$. We then denote its {\em translation
  axis} by
\[
\Ax_\ga=\{x\in X:d(x,\ga x)=\lambda(\ga)\}\,,
\]
and its repelling, attracting fixed points at infinity by $\ga_-,\ga_+
\in \partial_\infty X$ respectively. The element $\ga$ is {\em
  elliptic} if it has a fixed point in $X$ and {\em parabolic} if it
is neither loxodromic nor elliptic.

The discrete group $\Ga$ acts by conjugation on the nonempty set
$I_\Ga$ of involutions of $\Ga$.  In what follows, $I$ and $J$ will
denote two fixed $\Ga$-invariant nonempty subsets of $I_\Ga$, and will
be endowed with the left action by conjugation of $\Ga$.  For
instance, given an involution $\alpha$ of $\Ga$, the set of the
conjugates of $\alpha$ by the elements of $\Ga$ will be denoted by
$I_\alpha$.

Let $\alpha\in I_\Ga$. We denote by
\[
F_\alpha=\{x\in X:\alpha x =x\}
\]
the fixed point set of $\alpha$, which is a proper
nonempty\footnote{For every $x\in X$, since $\alpha$ is isometric and
has order $2$, the midpoint of $[x,\alpha x]$ is fixed by $\alpha$.}
closed convex\footnote{If $x,y\in F_\alpha$, then $\alpha([x,y])=
[x,y]$ as $\alpha$ is an isometry and $X$ is uniquely geodesic, hence
$[x,y]\subset F_\alpha$.}  subset of $X$. It is a totally geodesic
submanifold when $X$ is a manifold by for instance \cite[\S 2.80
  bis]{GalHulLaf90}, and the geometric realisation of a simplicial
subtree of $\XX$ when $X$ is a tree.

\brema\label{rem:remdefinvol} \hypertarget{remdefinvol1}{(1)} If
$X=\wt M$ is manifold with dimension $m\geq 1$, for $\alpha\in
\Isom(X)$ an involution, the dimension $k$ of the submanifold
$F_\alpha$ could be any element $k\in\llbracket 0, m-1
\rrbracket$,\footnote{In this paper, for all $a,b\in\ZZ$ with $a\leq
b$, we denote $\llbracket a,b\rrbracket=[a,b]\cap\ZZ$.}  as can be
seen when $\wt M$ is the real hyperbolic $m$-space $\HH^m_\RR$.  For
every $x\in F_\alpha$, the tangent space $T_x \wt M$ decomposes as an
orthogonal sum $T_x\wt M= T_xF_\alpha \oplus \nu_xF_\alpha$ of the
$k$-dimensional tangent space $T_xF_\alpha$ and the
$(m-k)$-dimensional normal space $\nu_x F_\alpha$ at $x$ to the fixed
point set $F_\alpha$. The tangent map $T_x\alpha$ acts by the block
matrix $\begin{psmallmatrix}\id&\ \ 0\\0& -\id\end{psmallmatrix}$ on
  this decomposition of $T_x\wt M$.  Thus, the involution $\alpha$
  reverses the orientation of $\wt M$ if and only if $m-k$ is odd.

\medskip
\hypertarget{remdefinvol2}{(2)} If $X=\wt M$ is a manifold, then the
action of any $\alpha\in I_\Ga$ on $\wt M$ is determined by $F_\alpha$:
\[
\forall\;\alpha,\beta\in I_\Ga,\qquad\text{if}\quad
F_\alpha=F_\beta\quad\text{then}\quad\alpha=\beta\,.
\]
Indeed, if $x\in F_\alpha$, then $\alpha x=x$, and if $x\in\wt M\ssm
F_\alpha$, if $p$ is the closest point to $x$ on $F_\alpha$, then
$\alpha x$ is the symmetric point of $x$ with respect to $p$ on the
(unique) geodesic line through $x$ and $p$. But this is no longer true
when $X$ is a tree, see the comment after Lemma
\ref{lem:cardequivmodI}, and Lemma \ref{lem:fixsettree}.
\erema

An element $\ga\in\Ga$ is {\it reversible} in $\Ga$ if it is
conjugated to its inverse.  It is {\it strongly reversible} in $\Ga$
if it is conjugated to its inverse by an involution in $\Ga$. We refer
to \cite[Sect.~2]{OFaSho15} for the basic ideas of reversibility, and
the rest of the cited book for an extensive survey. In particular, an
element of $\Ga$ is strongly reversible if and only if it is the
product of two involutions of $\Ga$, see \cite[Prop.~2.12]{OFaSho15}.

In this paper, we are
interested in strongly reversible loxodromic elements of $\Ga$. 
An element $\ga\in\Ga$ is {\it $\{I,J\}$-reversible} if there exists
an element $\alpha\in I$ such that we have $\ga\alpha =\alpha\ga^{-1}
\in J$. Such an element $\alpha\in I_{\Ga}$ is called a {\it
  $\ga$-reversing involution for $(I,J)$}.   We
denote by $\wt {\mathfrak R}_{I,J}$ the set of $\{I,J\}$-reversible
loxodromic elements of $\Ga$. Note that if $I'$ and $J'$ are
$\Ga$-invariant subsets of $I_\Ga$ such that $I\subset I'$ and
$J\subset J'$, then $\wt {\mathfrak R}_{I,J}\subset\wt {\mathfrak
  R}_{I',J'}$.

\brema\label{rem:defIrecip1} \hypertarget{defIrecip1_1}{(1)} Given a
$\ga$-reversing involution $\alpha$ for $(I,J)$, the set of
$\ga$-reversing involutions for $(I,J)$ is equal to $\{\alpha'\in
(\alpha Z_\Ga(\ga) )\cap I:\ga\alpha'\in J\}$. Since
$\beta\alpha=(\beta\alpha\beta^{-1})\beta$ for all $\alpha\in I$ and
$\beta\in J$, the element $\beta$ is a $\ga$-reversing involution for
$(J,I)$. Therefore being $\{I,J\}$-reversible or $\{J,I\}$-reversible is
equivalent, thus explaining the notation.  In particular, the set $\wt
{\mathfrak R}_{I,J}$ is invariant under taking inverses.

\hypertarget{defIrecip1_2}{(2)} The left action of $\Ga$ on itself by
conjugation preserves $\wt {\mathfrak R}_{I,J}$, since for every
$\delta \in\Ga$, an involution $\alpha$ is $\ga$-reversing for $(I,J)$
if and only if the element $\delta \alpha\delta^{-1}$, which belongs
to $I$, is $\delta\ga \delta^{-1}$-reversing for $(I,J)$, as
\[(\delta \alpha
\delta^{-1})(\delta \ga \delta^{-1})(\delta \alpha\delta^{-1})^{-1} =
(\delta \ga \delta^{-1}) ^{-1}
\quad\text{and}\quad
(\delta\ga \delta^{-1})(\delta
\alpha\delta^{-1}) = \delta (\ga\alpha) \delta^{-1}\in J\;.
\] 

\hypertarget{defIrecip1_3}{(3)} The set ${\mathfrak R}_{I,J}$ of
$\{I,J\}$-reversible loxodromic elements of $\Ga$ is invariant by
taking odd powers: Indeed, let $\ga\in \wt {\mathfrak R}_{I,J}$ and
let $\alpha$ be a $\ga$-reversing involution for $(I,J)$, so that
$\alpha\in I$, $\beta=\ga\alpha\in J$ and $\ga= \beta\alpha$. Since
$\alpha\beta= \alpha^{-1}\beta^{-1}= (\beta\alpha)^{-1}$ and
since $J$ is stable by conjugation, for every $k\in\NN\ssm\{0\}$, we
have
\[
\ga^{2k-1}=(\beta\alpha)^{2k-1}=((\beta\alpha)^{k-1}\beta
(\beta\alpha)^{-k+1})\;\alpha\in\wt {\mathfrak R}_{I,J}\;.
\]
On the other hand, since $I$ and $J$ are stable by conjugation, we have
\begin{align*}
  \ga^{2k}=(\beta\alpha)^{2k} & =(\beta\alpha\beta^{-1})
  \big((\alpha\beta)^{k-1}\alpha(\alpha\beta)^{-k+1}\big)\\
& =\big((\beta\alpha)^{k-1}\beta(\beta\alpha)^{-k+1}\big) (\alpha
\beta\alpha^{-1})\in\wt {\mathfrak R}_{I,I}\cap\wt {\mathfrak R}_{J,J}\,.
\end{align*}

\hypertarget{remdefIrecip4}{(4)} Let $\ga$ be a reversible loxodromic
element of $\Ga$.  If $g$ is an isometry such that $g\ga
g^{-1}=\ga^{-1}$, then $g$ preserves the translation axis of $\ga$ and
exchanges its two endpoints at infinity. In particular, the
restriction of $g$ to $\Ax_\ga$ is an orientation reversing isometry
of the geodesic line $\Ax_\ga$.  Hence, $g$ has a unique fixed point
$f_{g,\ga}$ on $\Ax_\ga$, so that $\{f_{g,\ga}\}= F_g\cap\Ax_\ga$.  In
particular, $g$ is elliptic. The two open subrays of $\Ax_\ga$ defined
by removing the point $f_{\alpha,\ga}$ are exchanged by $\alpha$. In
particular, the order of $\alpha$ is even.

When $X$ is a manifold, the two opposite unit tangent
vectors to $\Ax_\ga$ at $f_{\alpha, \ga}$, which are normal to
$F_\alpha$, are exchanged by $\alpha$. 

\medskip
\hypertarget{remdefIrecip2}{(5)} When $X=\wt M$ is a manifold with
dimension $2$ and when $\Ga$ is contained in the orientation
preserving isometry group of $\wt M$, as in \cite{Sarnak07,ErlSou24},
the fixed point sets of involutions are singletons, and every
reversible loxodromic element $\ga$ of $\Ga$ is strongly
reversible.\footnote{This result is mentioned without proof for the
hyperbolic plane in \cite[Sect.~6.2]{OFaSho15}.}  Indeed, if $g\in\Ga$
is an element such that $g\ga g^{-1}=\ga^{-1}$, then $g$ is elliptic
and preserves $\Ax_\ga$ by Remark (\hyperlink{remdefIrecip4}{4}).
Furthermore, $g$ cannot have a rotational component around $\Ax_\ga$
since $\wt M$ has dimension $2$ and $g$ preserves the orientation.

In general, an elliptic element $\alpha\in \Ga$ such that $\ga^{-1}=
\alpha\ga \alpha^{-1}$ can have an arbitrary even order. For instance,
for every $\lambda>1$, the Poincar\'e extension of the homography
$z\mapsto \lambda z$ to the upper halfspace model $\htr$ has a trivial
rotation factor around $\Ax_\ga$, which is the geodesic line with
points at infinity $0$ and $\infty$.  For every $k\in\NN\ssm\{0\}$,
the unique fixed point of the Poincar\'e extension of the mapping
$\alpha: z\mapsto e^{\frac{2\pi i}k}\;{\overline z}^{\,-1}$ is
$(0,1)\in\Ax_\ga$.  The isometry $\alpha$ is the composition of a
reflexion in the Euclidean unit sphere that exchanges the two
endpoints of $\Ax_\ga$, and a rotation of order $k$ around $\Ax_\ga$.
Hence $\ga^{-1}= \alpha\ga\alpha^{-1}$ and $\alpha$ has order $k$ if
$k$ is even and order $2k$ otherwise.
\erema

A subgroup $H$ of $\Ga$ isomorphic to the free product
$\ZZ/2\ZZ*\ZZ/2\ZZ$ of two copies of the cyclic group of order $2$ is
an {\it infinite dihedral subgroup} of $\Ga$.  Note that $H$ is
isomorphic to the semidirect product $\ZZ\rtimes\ZZ/2\ZZ$ and it has a
canonical morphism onto $\ZZ/2\ZZ$ (mapping any element of order $2$
to $1\in\ZZ/2\ZZ$), with kernel a canonical infinite cyclic (normal)
subgroup $Z_H$. If $\alpha\in I$, $\beta\in J$ and $H=\langle
\alpha\rangle *\langle \beta\rangle$, then $\ga=\beta\alpha $
generates $Z_H$, and $\ga$ is $\{I,J\}$-reversible in $\Ga$ with
$\alpha$ a $\ga$-reversing involution for $(I,J)$, since
\begin{equation}\label{eq:dihedralreciprocal}
  \ga^{-1}=\alpha\beta  =\alpha(\beta \alpha) \alpha
  =\alpha\ga\alpha^{-1}\,.
\end{equation}

Let $D^-$ and $D^+$ be two nonempty proper closed convex subsets of
$X$. A geodesic arc $\wt\rho: [0,T]\ra X$, where $T>0$, $\wt\rho(0)\in
D^-$ and $\wt\rho(T)\in D^+$ is a {\it common perpendicular of length
  $\lambda(\wt\rho)=T$ from $D^-$ to $D^+$} if its image is the unique
shortest geodesic segment from a point of $D^-$ to a point of
$D^+$. It exists if and only if the closures $\overline{D^-}$ and
$\overline{D^+}$ of $D^-$ and $D^+$ in $X\cup\partial_\infty X$ are
disjoint.  We refer to \cite{ParPau17ETDS}, the survey
\cite{ParPau16LMS}, or \cite[\S 2.5, \S 12.1 and \S 12.4]{BroParPau19}
for further details on common perpendiculars.

The following lemma connects strongly reversible loxodromic elements
of $\Ga$ with common perpendiculars of the fixed point sets of
involutions. This connection is the key to the proofs of the main
results in Sections \ref{sec:proofs} and \ref{sec:eqproof}.

\blemm\label{lem:prodinversion} Let $(\alpha,\beta)\in I_\Ga\times
I_\Ga$ and let $I$ and $J$ be $\Ga$-invariant nonempty subsets of
$I_\Ga$. The product $\ga=\beta\alpha$ is strongly reversible, and it
is $\{I,J\}$-reversible if furthermore $\alpha\in I$ and $\beta\in
J$. It is
\begin{enumerate}[label=\roman*]
\item[\hypertarget{itm:prod1}{(i)}]
  elliptic or the identity if and only if $F_\alpha\cap
    F_\beta\neq\emptyset$,
  \item[\hypertarget{itm:prod2}{(ii)}]   loxodromic if and only if
    $\overline{F_\alpha} \,\cap\,\overline{ F_\beta} =\emptyset$,
  \item[\hypertarget{itm:prod3}{(iii)}] parabolic if and only if
    $F_\alpha\cap F_\beta=\emptyset$ and $\overline{F_\alpha}\,\cap
    \,\overline{F_\beta}\neq\emptyset$.
\end{enumerate}

In cases (\hyperlink{itm:prod2}{ii}) and (\hyperlink{itm:prod3}{iii}),
the group $\langle\alpha,\beta\rangle$ generated by $\alpha$ and
$\beta$ is an infinite dihedral group and we have
\[
\langle\alpha,\beta\rangle=\langle\alpha\rangle*\langle\beta\rangle\;.
\]
If $\ga$ is loxodromic, the translation axis $\Ax_\ga$ of $\ga$ is the
union of the images by the elements of the group $\langle \alpha,
\beta\rangle$ of the common perpendicular segment $\wt\rho_{\alpha,
  \beta}$ from $F_{\alpha}$ to $F_{\beta}$, and the translation length
of $\ga$ satisfies
\begin{equation}\label{eq:length2length}
\lambda(\ga)=2\;d(F_\alpha,F_\beta)\,.
\end{equation}
\elemm

\begin{center}
  \begin{picture}(0,0)%
\includegraphics{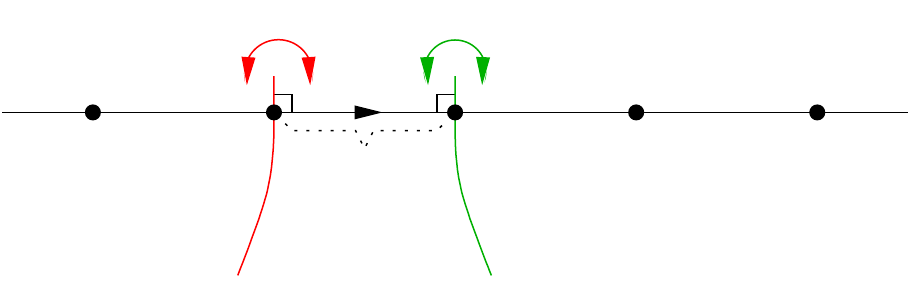}%
\end{picture}%
\setlength{\unitlength}{3812sp}%
\begingroup\makeatletter\ifx\SetFigFont\undefined%
\gdef\SetFigFont#1#2#3#4#5{%
  \reset@font\fontsize{#1}{#2pt}%
  \fontfamily{#3}\fontseries{#4}\fontshape{#5}%
  \selectfont}%
\fi\endgroup%
\begin{picture}(4524,1426)(439,-485)
\put(1666,254){\makebox(0,0)[lb]{\smash{{\SetFigFont{11}{13.2}{\rmdefault}{\mddefault}{\updefault}{\color[rgb]{0,0,0}$x$}%
}}}}
\put(2746,254){\makebox(0,0)[lb]{\smash{{\SetFigFont{11}{13.2}{\rmdefault}{\mddefault}{\updefault}{\color[rgb]{0,0,0}$y$}%
}}}}
\put(766,209){\makebox(0,0)[lb]{\smash{{\SetFigFont{11}{13.2}{\rmdefault}{\mddefault}{\updefault}{\color[rgb]{0,0,0}$\alpha y$}%
}}}}
\put(1756,794){\makebox(0,0)[lb]{\smash{{\SetFigFont{11}{13.2}{\rmdefault}{\mddefault}{\updefault}{\color[rgb]{1,0,0}$\alpha$}%
}}}}
\put(1711,-421){\makebox(0,0)[lb]{\smash{{\SetFigFont{11}{13.2}{\rmdefault}{\mddefault}{\updefault}{\color[rgb]{1,0,0}$F_\alpha$}%
}}}}
\put(2926,-421){\makebox(0,0)[lb]{\smash{{\SetFigFont{11}{13.2}{\rmdefault}{\mddefault}{\updefault}{\color[rgb]{0,.69,0}$F_\beta$}%
}}}}
\put(2656,794){\makebox(0,0)[lb]{\smash{{\SetFigFont{11}{13.2}{\rmdefault}{\mddefault}{\updefault}{\color[rgb]{0,.69,0}$\beta$}%
}}}}
\put(2071, 74){\makebox(0,0)[lb]{\smash{{\SetFigFont{11}{13.2}{\rmdefault}{\mddefault}{\updefault}{\color[rgb]{0,0,0}$\wt\rho_{\alpha,\beta}$}%
}}}}
\put(3601,209){\makebox(0,0)[lb]{\smash{{\SetFigFont{11}{13.2}{\rmdefault}{\mddefault}{\updefault}{\color[rgb]{0,0,0}$\beta x$}%
}}}}
\put(4501,209){\makebox(0,0)[lb]{\smash{{\SetFigFont{11}{13.2}{\rmdefault}{\mddefault}{\updefault}{\color[rgb]{0,0,0}$\beta\alpha y$}%
}}}}
\end{picture}%

\end{center}

\dem Let $\alpha\in I$ and $\beta\in J$ (with possibly
$I=J=I_\Ga$). The element $\alpha$ is $(I,J)$-reversing for
$\ga=\beta\alpha$, and hence $\ga$ is $\{I,J\}$-reversible, since
\begin{equation}\label{eq:prodinvrecloxo}
\alpha \ga\alpha^{-1}= \alpha (\beta \alpha ) \alpha^{-1}=\alpha
\beta=(\beta \alpha )^{-1}=\ga^{-1}\,.
\end{equation}

If $F_\alpha\cap F_\beta\neq\emptyset$, then $\ga=\beta\alpha$
pointwise fixes this intersection and $\ga$ is elliptic.

Assume that $F_\alpha\cap F_\beta=\emptyset$. If the intersection
$\overline{F_\alpha} \,\cap\,\overline{F_\beta}$ is nonempty, then by
convexity, this intersection is reduced to a point at infinity
$\xi\in\partial_\infty X$, which is fixed by $\alpha$ and $\beta$,
hence by $\ga$. If $X$ is a tree, then the two convex subsets
$F_\alpha$ and $F_\beta$ would meet, which has been excluded. Hence
$X$ is a manifold. Any horosphere centred at $\xi$ is invariant by
$\alpha$ and $\beta$, hence by $\ga=\beta\alpha$ and by $\ga^{-1}=
\alpha\beta$, so that $\ga$ is parabolic or elliptic.  If a point
$x\in X$ is fixed by $\beta\alpha$, then $\beta$ sends the segment
$[x,\alpha x]$ to the segment $[\beta x,x]$, hence it sends the
midpoint $m_\alpha$ of $[x,\alpha x]$ to the midpoint $m_\beta$ of
$[\beta x,x]$.  But $m_\alpha\in F_\alpha$ and $m_\beta\in F_\beta$.
Hence $m_\alpha=\beta^{-1}m_\beta=m_\beta\in F_\alpha \cap F_\beta$,
contradicting the disjointness of $F_\alpha$ and $F_\beta$. Thus $\ga$
is not elliptic, therefore it is parabolic.  In particular, $\ga$ has
infinite order, and $\langle \alpha, \beta\rangle$ is indeed
isomorphic to $\langle\alpha \rangle *\langle\beta\rangle$.

If $\overline{F_\alpha} \,\cap\,\overline{ F_\beta} =\emptyset$, let
$\wt\rho_{\alpha,\beta}=[x,y]$ be the common perpendicular from
$F_\alpha$ to $F_\beta$ (with $x\in F_\alpha$ and $y\in F_\beta$, see
the above picture). We claim that $\alpha [y,x]\cup [x,y]$ is a
geodesic segment. When $X$ is a manifold, this follows from the
properties of the action of $T_x\alpha$ on $\nu_xF_\alpha$ described
in Remark \ref{rem:defIrecip1} (\hyperlink{remdefIrecip2}{5}).  When
$X$ is a tree, the interior of $[x,y]$ would otherwise contain a fixed
point of $\alpha$. Similarly, $[\alpha y,y]\cup\beta [y,\alpha y]$ is
a geodesic segment, and $\ga=\beta\alpha$ maps its first half $[\alpha
  y,y]$ to its second half $[\beta y,\beta \alpha y]$ in an
orientation preserving way. This implies that $\ga$ is loxodromic,
that its translation axis contains $[\alpha y,y]$, hence contains
$\wt\rho_{\alpha,\beta}$, and that
\[
\lambda(\ga)=d(\alpha y,y)=2 \,d(x,y)=2\,d(F_\alpha,F_\beta)\,.
\]
Taking the image of the common perpendicular segment $\wt
\rho_{\alpha,\beta}$ by $\alpha$, and then the images of $\alpha\wt
\rho_{\alpha,\beta}\cup\wt \rho_{\alpha,\beta}$ by the powers of
$\ga$, we cover the whole translation axis $\Ax_\ga$.  The fact that
$\langle\alpha,\beta \rangle$ is isomorphic to $\langle\alpha\rangle
*\langle\beta\rangle$ follows as in the parabolic case.
\cqfd

\medskip
We say that an infinite dihedral subgroup $H$ of $\Ga$ is {\it of
  loxodromic type} if $Z_H$ is generated by a loxodromic element and
that it is {\it of parabolic type} if $Z_H$ is generated by a
parabolic element. Note that when $X$ is a manifold with dimension $2$
and when $\Ga$ preserves an orientation, then all infinite dihedral
subgroups of $\Ga$ are of loxodromic type.  In case
(\hyperlink{itm:prod2}{ii}) of Lemma \ref{lem:prodinversion}, the
group $\langle\alpha,\beta\rangle$ generated by $\alpha$ and $\beta$
is an infinite dihedral group of loxodromic type.  In case
(\hyperlink{itm:prod3}{iii}), it is an infinite dihedral group of
parabolic type.

An {\it $\{I,J\}$-dihedral subgroup} of $\Ga$ is an infinite dihedral
subgroup $\langle \alpha\rangle *\langle \beta\rangle$ of loxodromic
type with $\alpha\in I$ and $\beta\in J$. Then $\ga= \beta\alpha$ is
$\{I,J\}$-reversible by Equation \eqref{eq:dihedralreciprocal} and
since $\ga\alpha=\beta\in J$.  We denote by $\wt {\mathfrak D}_{I,J}$
the set of $\{I,J\}$-dihedral subgroups of $\Ga$. It is invariant
under conjugation by elements of $\Ga$.

\bexer A classical example considered by Sarnak \cite{Sarnak07} is
$\Ga=\PSL_2(\ZZ)$ acting by homographies on the upper halfplane model
of the real hyperbolic plane $\wt M=\hdr$, for which $I_\Ga$ is the
set $I_\alpha$ of conjugates in $\Ga$ of $\alpha=\pm
\begin{psmallmatrix}\ \,0&1\\-1&0\end{psmallmatrix}\in\PSL_2(\ZZ)$.
As said in the Introduction, Sarnak calls the reversible loxodromic
elements of $\PSL_2(\ZZ)$ reciprocal. They are strongly reversible by
Remark \ref{rem:defIrecip1} (\hyperlink{remdefIrecip2}{5}).

Let $X = \hdr$ and let $\Ga$ be the extended modular group, generated
by the modular group $\PSL_2(\ZZ)$ acting by homographies on $\hdr$
and by the hyperbolic reflexion $\alpha:z\mapsto -\overline{z}$ fixing
the vertical geodesic line $F_\alpha$ with points at infinity
$0,\infty\in\partial_\infty\hdr=\RR\cup\{\infty\}$. Let $\beta:z
\mapsto \frac{-2\,\overline{z}+1}{-3\,\overline{z}+2}\in \Ga$, which
is the hyperbolic reflexion fixing the geodesic line $F_{\beta}$ with
points at infinity $\frac{1}{3}$ and $1$.  The composition
$\ga=\pm\begin{psmallmatrix}2&1\\3&2\end{psmallmatrix} =
\beta\alpha\in\Ga$ is an {\em ambiguous}\footnote{See \cite{Sarnak07},
\cite[\S 6]{ParPau25a}.} loxodromic element of $\PSL_2(\ZZ)$.  It is
$\{I_\alpha,I_\beta\}$-reversible since $\alpha\in I_\alpha$
conjugates $\ga$ to its inverse $\ga^{-1}=\alpha\beta$ and
$\beta=\ga\alpha\in I_\beta$.  But $\ga$ is not $\{I_\alpha,
I_\alpha\}$-reversible, since $\beta=\ga\alpha$ is not conjugated to
$\alpha$ in $\Ga$, hence does not belong to $I_{\alpha}$. Note that
$\ga^2= (\beta \alpha \beta^{-1})\alpha= \beta (\alpha
\beta\alpha^{-1})$ is $\{I_\alpha, I_\alpha\}$-reversible and
$\{I_\beta,I_\beta\}$-reversible as said in Remark
\ref{rem:defIrecip1} (\hyperlink{defIrecip1_3}{3}).  The translation
axis of $\ga$ is the geodesic line in $\hdr$ containing the common
perpendicular from $F_\alpha$ to $F_\beta$ by Lemma
\ref{lem:prodinversion}.  See the end of Example
\ref{ex:extendedmodular} and \cite[\S 6]{ParPau25a} for further
details on this example.
\eexer

\medskip
We refer to \cite[\S 7.2]{BroParPau19} for the following definitions
concerning equivariant families of (convex) subsets of $X$.  Let $I'$
be an index set endowed with a left action of $\Ga$. A family
$\D=(D_i)_{i\in I'}$ of subsets of $X$ indexed by $I'$ is {\em
  $\Ga$-equivariant} if $\ga D_i=D_{\ga i}$ for all $\ga\in\Ga$ and
$i\in I'$. We denote by $\Sim[1.3]\D$ the equivalence relation on $I'$
defined by $i\Sim[1.3]\D j$ if and only if $D_i=D_j$ and there exists
$\ga\in\Ga$ such that $j=\ga i$.  We say that $\D$ is {\it locally
  finite} if for every compact subset $K$ in $X$, the quotient set
$\{i\in I': D_i\cap K\neq\emptyset\}/\!_{\Sim\D}$ is finite.

The family $\F_I=(F_\alpha)_{\alpha\in I}$ is $\Ga$-equivariant since
$\ga F_\alpha=F_{\ga\alpha\ga^{-1}}$ for all $\ga\in\Ga$ and
$\alpha\in I$, and locally finite by the discreteness of $\Ga$. We
simplify the notation $\Sim[1.4]{\F_I}$ as $\sim_I$, with 
\begin{equation}\label{def:simi}
\text{$\alpha\sim_I
\beta$ if and only if $F_\alpha=F_\beta$ and there exists $\ga\in\Ga$
such that $\alpha=\ga \beta\ga^{-1}$.} 
\end{equation}
When $X$ is a manifold, by Remark \ref{rem:remdefinvol}
(\hyperlink{remdefinvol2}{2}), we have $\alpha\sim_I \beta$ if and
only if $\alpha= \beta$.

The equivalence relation $\sim_I$ is $\Ga$-equivariant: for all
$\alpha,\beta\in I$ and $\ga\in\Ga$, we have $\ga\alpha\ga^{-1}
\sim_I\ga\beta\ga^{-1}$ if and only if $\alpha \sim_I\beta$. We
henceforth endow the quotient set $I/\!_{\sim_I}$ with the quotient
left action of the group $\Ga$.

For every $\alpha\in I_\Ga$, let $\Ga_{F_\alpha}=\stab_\Ga (F_\alpha)$
be the (global) stabiliser of $F_\alpha$ in $\Ga$. Let $Z_\Ga(\alpha)$
be the centraliser of $\alpha$ in $\Ga$, which is contained in
$\Ga_{F_\alpha}$.

\blemm\label{lem:decompreciloxo} For every $\ga\in\wt{\mathfrak
  R}_{I,J}$, the set of pairs $(\alpha,\beta)\in I\times J$ such that
$\ga= \beta \alpha$ is nonempty and invariant under the diagonal left
action by conjugation on the product $I\times J$ of the centraliser
$Z_\Ga(\ga)$ of $\ga$, and has finite quotient by this action.
\elemm

\dem Let $\ga\in \wt{\mathfrak R}_{I,J}$. By definition, there exists
$\alpha\in I$ such that $\beta=\ga\alpha\in J$, hence
$\ga=\beta\alpha^{-1}=\beta\alpha$. For every $\delta\in Z_\Ga(\ga)$,
we have
\[\ga=\delta\ga\delta^{-1}=
\delta\beta \alpha \delta^{-1}=(\delta\beta\delta^{-1}) (\delta\alpha
\delta^{-1})\,,
\]
which proves the first claim since $I$ and $J$ are invariant under
conjugation.

To prove the second claim, let us fix $x\in\Ax_\ga$.  Since
$Z_\Ga(\ga)$ contains $\ga^\ZZ$ which acts on $\Ax_\ga$ with
fundamental domain the relatively compact geodesic segment $[x,\ga
x[\,$, up to conjugating $\alpha$ by an element of $Z_\Ga(\ga)$, we
may assume that the unique fixed point $f_{\alpha, \ga}$ of $\alpha$
on $\Ax_\ga$ belongs to $[x,\ga x]$. The claim follows since the
family $\F_I$ is locally finite.
\cqfd

\blemm\label{lem:cardequivmodI} The group $Z_\Ga(\alpha)$ has finite
index in $\Ga_{F_\alpha}$.  If $\alpha\in I$, the map from $\Ga_{F_\alpha}/
Z_\Ga(\alpha)$ to the equivalence class in $I/\!_{\sim_I}$ of $\alpha$
 defined by $\ga Z_\Ga(\alpha)\mapsto \ga \alpha\ga^{-1}$
is a bijection.
\elemm

\dem By the discreteness of $\Ga$, the point stabilizers of $\Ga$ in
$X$ are finite, hence the equivalence classes of $\sim_I$ are finite.
Therefore the first claim follows from the second one.  The map from
$\Ga/Z_\Ga(\alpha)$ to the conjugacy class $[\alpha] =\{\ga \alpha
\ga^{-1}:\ga\in\Ga\}$ of $\alpha$ in $\Ga$ defined by $\ga
Z_\Ga(\alpha)\mapsto \ga \alpha\ga^{-1}$ is a bijection. For every
$\ga\in\Ga$, we have $F_{\ga\alpha\ga^{-1}}= F_\alpha$ if and only if
$\ga\in \Ga_{F_\alpha}$. The result follows by the definition of
$\sim_I$.
\cqfd

\medskip
Note that when $X$ is a tree, we can have $|\Ga_{F_\alpha}/
Z_\Ga(\alpha)| >1$, see Lemma \ref{lem:fixsettree}.  However, when $X$
is a manifold, we have $|\Ga_{F_\alpha}/ Z_\Ga(\alpha)|=1$ and the
equivalence classes for $\sim_I$ are singletons.

\section{Multisets of strongly reversible loxodromic elements}
\label{sec:multi}

A {\em set with multiplicity} or a {\em (positive real valued)
  multiset} is a pair $A=(\underline A,\omega)$, where $\underline A$
is a set, called the {\em underlying set} of $A$, and $\omega:
\underline A\to\interval[open]0{+\infty}$ is a positive function,
called the {\em multiplicity} of $A$. By an {\it element of} $A$, we
mean an element of $\underline A$. If $E=(\underline E,\omega)$ and
$E'= (\underline E',\omega')$ are multisets, a {\it multiset
  bijection} $f:E\to E'$ is a surjective map $f:\underline E\ra
\underline E'$ such that for every $x'\in E'$, we have $$\omega'(x')=
\sum_{x\in f^{-1} (x')} \omega(x)\,.$$ If $E=(\underline E,\omega)$ is
a multiset and $\underline A\subset\underline E$, the multiset
$A=(\underline A,\omega|_A)$ is a {\em restriction} of $E$.

As a preparation for the proofs in Sections \ref{sec:proofs} and
\ref{sec:eqproof}, we define various multisets associated with
strongly reversible loxodromic elements of $\Ga$.  At the end of this
Section, we prove that there are multiset bijections between these
multisets.

Let 
\[
\wt{IJ}=\big\{(\alpha,\beta)\in I\times J :
\overline{F_{\alpha}}\,\cap\,\overline{F_{\beta}}\,=\emptyset\big\}\,.
\]
Note that $(\alpha,\beta)\in\wt{IJ}$ if and only if $(\beta,\alpha)
\in\wt{JI}$.  Let $\sim_{\!IJ}$ be the diagonal equivalence relation
on $I\times J$, defined by $(\alpha,\beta)\sim_{\!IJ} (\alpha',
\beta')$ if and only if $\alpha\sim_I\alpha'$ and $\beta \sim_{J}
\beta'$. The condition $\overline{F_{\alpha}}\,\cap\,
\overline{F_{\beta}}\, =\emptyset$ for $(\alpha, \beta)\in I\times J$
is constant on the equivalence classes modulo $\sim_{\!IJ}$ and on the
orbits of the diagonal left action of the group $\Ga$ on $I \times J$.
Hence the subset $\wt{IJ}$ of $I\times J$ is saturated by
$\sim_{\!IJ}$ and $\Ga$-invariant.

We define the multiset $IJ$ with underlying set $\underline{IJ}=
\Ga\bs\wt{IJ}$ and multiplicity given by
\begin{equation}\label{eq:multcomperpIIGa}
\forall\;\Ga(\alpha,\beta)\in IJ,\quad \mult(\Ga(\alpha,\beta))=
\frac{1}{\card(\Ga_{F_\alpha}\cap\Ga_{F_\beta})\;
    |\Ga_{F_\alpha}/Z_\Ga(\alpha)|\; |\Ga_{F_\beta}/Z_\Ga(\beta)|}\,.
\end{equation}
The right hand side of the above formula for $(\alpha, \beta)\in \wt
{IJ}$ is constant on the equivalence classes modulo $\sim_{\!IJ}$ and
on the orbits of $\Ga$ on $\wt{IJ}$, hence $\mult(\Ga(\alpha,\beta))$
is well defined.

For every $T>0$, we define the multiset restriction $IJ(T)$ of $IJ$ by
setting
\begin{equation}\label{eq:defiIIGaT}
\underline{IJ(T)}=
\big\{\Ga(\alpha,\beta)\in IJ:d(F_\alpha,F_\beta)\leq T\}\,.
\end{equation}

\medskip
We denote by ${\mathfrak R}_{I,J}$ the multiset, whose underlying set
is the set $\underline{{\mathfrak R}_{I,J}}=\Ga\bs\wt{\mathfrak R}_{I,
  J}$ of conjugacy classes $[\ga]$ of $\{I,J\}$-reversible loxodromic
elements $\ga$ of $\Ga$, and where the multiplicity of $[\ga]\in
{\mathfrak R}_{I,J}$ is defined by
\begin{align}\label{eq:multreclox}
  &\mult([\ga])=
  \sum_{\Ga(\alpha,\,\beta)\,\in\,IJ\;:\;\beta\alpha\,\in\,[\ga]}\;
  \mult(\Ga(\alpha,\beta))\,.
\end{align}
We have ${\mathfrak R}_{I,J}={\mathfrak R}_{J,I}$. For every $T >0$,
we denote by ${\mathfrak R}_{I,J}(T)$ the multiset restriction of
${\mathfrak R}_{I,J}$ to its elements with translation length at most
$T$.

\brema\label{rem:multiset} \hypertarget{remmultiset1}{(1)} Equation
\eqref{eq:multreclox} is independent of the choices of representatives
$(\alpha, \beta)\in \wt{IJ}$ in their $\Ga$-orbit, since $(\ga'\beta
{\ga'}^{-1})(\ga'\alpha{\ga'}^{-1})=\ga'(\beta\alpha){\ga'}^{-1}$ for
all $\ga'\in\Ga$.  The sum in Equation \eqref{eq:multreclox} has a
finite index set by the last assertion of Lemma
\ref{lem:decompreciloxo}. This sum is positive by the first
assertion of Lemma \ref{lem:decompreciloxo}.

\smallskip
\hypertarget{remmultiset1}{(2)} The multiplicity $\mult([\ga])$
depends only on the conjugacy class $[\ga]= [\ga^{-1}]$ of the
$\{I,J\}$-reversible loxodromic element $\ga$ in $\Ga$.

\smallskip
\hypertarget{remmultiset1}{(3)} Assume in this remark that $X=\wt M$
is a manifold and that the elliptic elements of $\Ga$ only have
isolated fixed points. Equivalently, no element of $\Ga\ssm\{\id\}$
fixes a nontrivial geodesic segment in $X$.  Thus, for all
$\alpha,\beta\in I_\Ga$ with $\alpha\ne\beta$, every element of
$\Ga_{F_\alpha}\cap \Ga_{F_\beta}$ is trivial, because it fixes
pointwise the (nontrivial) common perpendicular between $F_\alpha$ and
$F_\beta$. We have $|\Ga_{F_\alpha}/Z_\Ga(\alpha)|=|\Ga_{F_\beta}/
Z_\Ga(\beta)| =1$ as seen at the end of Section
\ref{sec:generalities}.  Hence we have
\[
\mult([\ga])=\card\big\{\Ga(\alpha,\,\beta)\,\in\,IJ
\;:\;[\beta\alpha]=[\ga]\big\}\,.
\]
\erema

Here is a particular case when the computation of the multiplicity of
the conjugacy class of a strongly reversible loxodromic is easy.  An
element $\ga\in \Ga$ is {\em primitive} if it is not a proper power in
$\Ga$.  An element $\ga\in \wt {\mathfrak R}_{I,J}$ is {\it
  $\{I,J\}$-primitive} if it is not a proper power of an element of
$\wt {\mathfrak R}_{I,J}$. This property being invariant under
conjugation by $\Ga$, we also say that the conjugacy class $[\ga]\in
{\mathfrak R}_{I,J}$ of $\ga$ is {\it $\{I,J\}$-primitive}. We denote
the set of conjugacy classes of $\{I,J\}$-primitive elements in $\Ga$
by ${\mathfrak R}_{I,\,J,\,\rm{prim}}$.

\blemm\label{lem:computmultprim} Assume that $\wt M$ has dimension
$2$, that $\Ga$ preserves an orientation of $\wt M$, and that
$I=J=I_\Ga$. If $\ga$ is a primitive strongly reversible loxodromic
element of $\Ga$, then
\[
  \mult([\ga])=2\,.
\]
\elemm

\dem Let $(\alpha,\,\beta)$ and $(\alpha',\,\beta')$ be two elements
of $\wt{IJ}$ such that $[\beta\alpha]=[\beta'\alpha']=[\ga]$.  Up to
replacing them by other elements in their $\Ga$-orbits in $\wt{IJ}$,
we may assume that $\beta\alpha=\beta'\alpha'=\ga$. Let $x_\alpha,
x_\beta, x_{\alpha'}, x_{\beta'}$ be the (isolated) fixed point of
$\alpha,\beta,\alpha', \beta'$ respectively. They belong to the
translation axis $\Ax_\ga$ of $\ga$. Since $[x_\alpha,\ga x_\alpha[$
    is a fundamental domain for the action of $\ga^\ZZ$ on $\Ax_\ga$,
    up to diagonally conjugating $(\alpha',\beta')$ by an element of
    $Z_\Ga(\ga)$, we may assume that $x_{\alpha'}\in [x_\alpha, \ga
      x_\alpha[\,$.
\begin{center}
\begin{picture}(0,0)%
\includegraphics{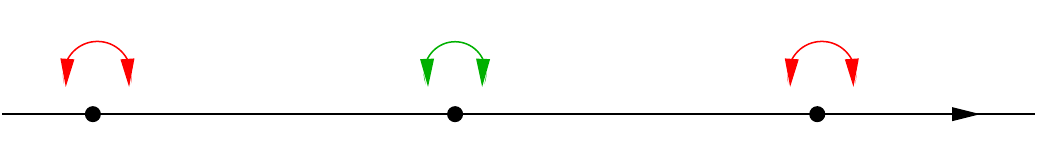}%
\end{picture}%
\setlength{\unitlength}{3812sp}%
\begingroup\makeatletter\ifx\SetFigFont\undefined%
\gdef\SetFigFont#1#2#3#4#5{%
  \reset@font\fontsize{#1}{#2pt}%
  \fontfamily{#3}\fontseries{#4}\fontshape{#5}%
  \selectfont}%
\fi\endgroup%
\begin{picture}(5247,808)(439,145)
\put(2656,794){\makebox(0,0)[lb]{\smash{{\SetFigFont{11}{13.2}{\rmdefault}{\mddefault}{\updefault}{\color[rgb]{0,.69,0}$\beta$}%
}}}}
\put(766,209){\makebox(0,0)[lb]{\smash{{\SetFigFont{11}{13.2}{\rmdefault}{\mddefault}{\updefault}{\color[rgb]{0,0,0}$x_\alpha$}%
}}}}
\put(856,794){\makebox(0,0)[lb]{\smash{{\SetFigFont{11}{13.2}{\rmdefault}{\mddefault}{\updefault}{\color[rgb]{1,0,0}$\alpha$}%
}}}}
\put(2566,209){\makebox(0,0)[lb]{\smash{{\SetFigFont{11}{13.2}{\rmdefault}{\mddefault}{\updefault}{\color[rgb]{0,0,0}$x_\beta$}%
}}}}
\put(4366,209){\makebox(0,0)[lb]{\smash{{\SetFigFont{11}{13.2}{\rmdefault}{\mddefault}{\updefault}{\color[rgb]{0,0,0}$\ga x_\alpha$}%
}}}}
\put(4366,794){\makebox(0,0)[lb]{\smash{{\SetFigFont{11}{13.2}{\rmdefault}{\mddefault}{\updefault}{\color[rgb]{1,0,0}$\beta\alpha\beta^{-1}$}%
}}}}
\put(5671,344){\makebox(0,0)[lb]{\smash{{\SetFigFont{11}{13.2}{\rmdefault}{\mddefault}{\updefault}{\color[rgb]{0,0,0}$\Ax_\ga$}%
}}}}
\end{picture}%

\end{center}
If $x_{\alpha'}=x_\alpha$, then since $\wt M$ is a manifold, we have
$\alpha'=\alpha$ by Remark \ref{rem:remdefinvol}
(\hyperlink{remdefinvol2}{2}).  Hence, $\beta'=\beta$ and
$\Ga(\alpha,\,\beta)= \Ga(\alpha',\,\beta')$. Otherwise, we have the
following three cases to consider.

If $x_{\alpha'}\in\;]x_\alpha, x_\beta[\,$, then the element $\alpha'
\alpha \in \Ga$ is loxodromic, with same translation axis as $\ga$,
but with translation length $\lambda(\alpha\alpha') = 2\,d(x_\alpha,
x_{\alpha'})<2\,d(x_\alpha,x_\beta)=\lambda(\ga)$, contradicting the
fact that $\ga$ is primitive.

If $x_{\alpha'}\in \;]x_\beta,\ga x_\alpha[\,$, then similarly
considering the element $\alpha'\beta$ contradicts the fact that $\ga$
is primitive.

Hence $x_{\alpha'}=x_\beta$. Therefore $\alpha'=\beta$, and
$x_{\beta'}=\ga x_\alpha=\beta x_\alpha=x_{\beta\alpha\beta^{-1}}$, so
that we have $\beta'=\beta\alpha \beta^{-1}$, and there is only one
possibility that $\Ga(\alpha',\, \beta')= \Ga(\beta,\alpha)$.  Note
that $\Ga(\alpha,\,\beta)\neq \Ga(\beta,\,\alpha)$: Otherwise there
exists $\delta\in\Ga$ so that $\beta=\delta\alpha\delta^{-1}$ and
$\alpha= \delta\beta\delta^{-1}$.  But then $\delta x_\alpha=x_\beta$
and $\delta x_\beta=x_\alpha$, so that $\delta$ fixes the midpoint of
the geodesic segment $[x_\alpha, x_\beta]$. Considering the element
$\delta\alpha$ again contradicts the fact that $\ga$ is primitive.
This proves the result.
\cqfd

\medskip
A {\it common perpendicular $\rho$ of type $(I,J)$} is the $\Ga$-orbit
of the common perpendicular $\wt \rho_{\alpha ,\beta}$ from
$F_{\alpha}$ to $F_{\beta }$, where $(\alpha, \beta)$ ranges over
$\wt{IJ}$. Note that
\begin{equation}\label{eq:saturatalphabeta}
\forall\,(\alpha,\beta),(\alpha',\beta')\in \wt{IJ},\qquad
\text{if}\quad (\alpha,\beta)\sim_{\!IJ}(\alpha',\beta')
\quad\text{then}\quad\wt\rho_{\alpha',\beta'} =\wt\rho_{\alpha,\beta}
\end{equation}
\begin{equation}\label{eq:equivarthoalphabeta}
 \text{and}\qquad
\forall\,(\alpha,\beta)\in \wt{IJ},\quad\forall\;\ga\in\Ga,\qquad
\ga \;\wt\rho_{\alpha,\beta} =\wt\rho_{\ga\alpha\ga^{-1},\ga\beta\ga^{-1}} \,.
\end{equation}
We then denote $\rho=\rho_{\alpha ,\beta}=\Ga\,\wt\rho_{\alpha,\beta}
\in\Ga\bs\gengeod X$. Remark that such a pair $(\alpha, \beta)$ is not
necessarily unique. The {\it length $\lambda (\rho)$} of $\rho$ is the
length of any such $\wt\rho_{\alpha, \beta}$.  The {\it multiplicity}
of $\rho$ is
\begin{equation}\label{eq:multcomperp}
 \mult(\rho) =\sum_{\Ga(\alpha,\,\beta)\,\in\, IJ, \;\;\rho=
   \rho_{\alpha,\beta}} \frac{1}{\card(\Ga_{F_\alpha}\cap\Ga_{F_\beta})\;
   |\Ga_{F_\alpha}/Z_\Ga(\alpha)|\; |\Ga_{F_\beta}/Z_\Ga(\beta)|}\,.
\end{equation}
This sum does not depend on the choice of a representative $(\alpha,
\beta)$ in each orbit of $\Ga$ in $\wt{IJ}$, and is finite by the
local finiteness of the families $\F_I$ and $\F_J$. Recall that if
$(\alpha,\beta)\in \wt{IJ}$, then $(\beta,\alpha)\in \wt{JI}$, and
note that
\begin{equation}\label{eq:multcomperpsym}
\mult(\rho_{\alpha,\beta})=\mult(\rho_{\beta,\alpha})\,.
\end{equation}

We denote by $\Perp(\F_I,\F_J)$ the multiset of the common
perpendiculars of type $(I,J)$. For every $T >0$, we denote by $\Perp
(\F_I,\F_J,T)$ the multiset restriction of $\Perp(\F_I,\F_J)$ to its
elements with length at most $T$ and by $\underline{\Perp}
(\F_I,\F_J,T)$ its underlying set.

\bprop\label{prop:bijcomperprecloxo} There exist two multiset
bijections $\Theta_1: IJ\to\Perp(\F_I,\F_J)$ and \mbox{$\Theta_2:
IJ\to{\mathfrak R}_{I,J}$} such that for every $(\alpha,\beta)\in
\wt{IJ}$, we have
\[
\Theta_1(\Ga (\alpha,\beta)) = \rho_{\alpha,\beta}
\quad\text{and}\quad\Theta_2(\Ga(\alpha,\beta))=[\beta\alpha]\,.
\]
Furthermore, for every $T>0$, the maps $\Theta_1$ and $\Theta_2$ send
$IJ(T)$ to $\Perp(\F_I,\F_J,T)$ and ${\mathfrak R}_{I,J}(2T)$ respectively.
\eprop

\dem The map $\Theta_1: \underline{IJ}=\Ga\bs \wt{IJ}\to
\underline{\Perp}(\F_I,\F_J)$ defined by $\Ga(\alpha,\beta)\mapsto
\rho_{\alpha,\beta}=\Ga\wt \rho_{\alpha,\beta}$ is well-defined by
Equation \eqref{eq:equivarthoalphabeta} and is surjective by the
definition of a common perpendicular of type $(I,J)$. By Equations
\eqref{eq:multcomperp} and \eqref{eq:multcomperpIIGa}, for every
element $\rho\in\Perp(\F_I,\F_J)$, we have $\mult(\rho) =
\sum_{\Ga(\alpha,\,\beta)\,\in\,\Theta_1^{-1}(\rho)}\mult(\alpha,\beta)$.
Hence $\Theta_1$ is a multiset bijection from $IJ$ to
$\Perp(\F_I,\F_J)$.

The map $\wt\Theta_2:\wt{IJ}\to\wt {\mathfrak R}_{I,J}$ defined by
\[
\wt\Theta_2:(\alpha ,\beta )\mapsto\beta \alpha 
\]
indeed has values in the set $\wt {\mathfrak R}_{I,J}$ of
$\{I,J\}$-reversible loxodromic elements of $\Ga$ by Lemma
\ref{lem:prodinversion}.  The map $\wt\Theta_2$ is equivariant for the
left actions of $\Ga$ by diagonal conjugation on $\wt{IJ}$ and by
conjugation on $\wt {\mathfrak R}_{I,J}$. The map $\wt\Theta_2$ is
surjective by Lemma \ref{lem:decompreciloxo}. Hence $\wt\Theta_2$
induces a surjective quotient map $\Theta_2: \underline{IJ} \ra
\underline{{\mathfrak R}_{I,J}}$ defined by $\Ga (\alpha,\beta)\mapsto
[\beta\alpha]$. By Equation \eqref{eq:multreclox}, we have
$\mult([\ga])= \sum_{\Ga(\alpha,\,\beta)\,\in\,(\Theta_2)^{-1}([\ga])}
\mult\Ga(\alpha, \beta)$. Hence $\Theta_2$ is a multiset bijection
from $IJ$ to ${\mathfrak R}_{I,J}$.

For all $(\alpha,\beta)\in \wt{IJ}$ and $T>0$, the common
perpendicular $\wt\rho_{\alpha ,\beta }$ from $F_{\alpha}$ to
$F_{\beta}$ has length $\lambda(\rho_{\alpha ,\beta })=d(F_\alpha,
F_\beta)$ at most $T$ if and only if the $\{I,J\}$-reversible
loxodromic element $\beta\alpha$ has translation length at most $2T$,
by Equation \eqref{eq:length2length}.  This proves the last claim of
Proposition \ref{prop:bijcomperprecloxo}.  \cqfd

\medskip
For every $\{I,J\}$-dihedral subgroup $H$ of $\Ga$, if $\ga$ generates
the infinite cyclic group $Z_H$, we define the multiplicity of $H$ as
\begin{equation}\label{eq:defimultdihedral}
\mult(H)=\mult([\ga])\,.
\end{equation}
This multiplicity does not depend on the choice of one of the two
generators $\ga$ of $Z_H$ since $[\ga]=[\ga^{-1}]$ and it is constant
on the conjugacy class of $H$.  We denote by ${\mathfrak D}_{I,J}$ the
multiset of conjugacy classes in $\Ga$ of $\{I,J\}$-dihedral subgroups
of $\Ga$.

\bprop\label{prop:bijdihedralrecloxo} The map $\wt\Theta_3$ from $\wt
{\mathfrak D}_{I,J}$ to ${\mathfrak R}_{I,J}$ which associates to
$H\in\wt {\mathfrak D}_{I,J}$ the conjugacy class of a generator of
$Z_H$ induces a multiset bijection $\Theta_3$ from the multiset
${\mathfrak D}_{I,J}$ to the multiset ${\mathfrak R}_{I,J}$.
\eprop

\dem We have already seen that for every $H\in\wt {\mathfrak
  D}_{I,J}$, any generator of $Z_H$ is $\{I,J\}$-reversible and hence,
so is its inverse.  Thus the map $\wt\Theta_3$ is well defined. It is
surjective by Lemmas \ref{lem:decompreciloxo} and
\ref{lem:prodinversion}.  It is invariant by conjugation since
$Z_{\ga'H{\ga'}^{-1}}= \ga'Z_H{\ga'}^{-1}$ for every $H\in\wt
{\mathfrak D}_{I,J}$. The definition of the multiplicities in Equation
\eqref{eq:defimultdihedral} proves that the induced map between the
sets with multiplicities ${\mathfrak D}_{I,J}$ and ${\mathfrak
  R}_{I,J}$ is a multiset bijection.
\cqfd

\medskip 
The map $\Theta_3$  restricts to a multiset bijection between the
multiset  of conjugacy classes of maximal\footnote{for
inclusion amongst $\{I,J\}$-dihedral subgroups of $\Ga$ }
$\{I,J\}$-dihedral subgroup of loxodromic type and the multiset
${\mathfrak R}_{I,\,J,\,{\rm prim}}$ of conjugacy classes of
$\{I,J\}$-primitive $\{I,J\}$-reversible loxodromic elements of $\Ga$.

\section{Dynamics of the geodesic flow}
\label{sec:dyn}

This whole section consists of background information from
\cite[Chap.~3, 6, 7]{PauPolSha15} and \cite[Chap.~3, 4]{BroParPau19},
to which we refer for proofs, details and complements, in particular
for the definitions of the measures generalizing the ones given in the
introduction.

\medskip
\noindent{\bf Geodesics and the geodesic flow.} Let $(X,\Ga)$ be as in
the beginning of Section \ref{sec:generalities}. Let $\G X$ be the
metric space of all geodesic lines $\ell: \RR\ra X$ in $X$, such that,
when $X$ is a tree, $\ell(0)$ is a vertex of $\XX$. The distance
between two elements $\ell$ and $\ell'$ in $\G X$ is
\begin{equation}\label{eq:ggdist}
d(\ell,\ell')= \int_{-\infty}^{\infty}
d(\ell(t),\ell'(t))\;e^{-2\,|t|}\;dt\,.
\end{equation}
When $X=\wt M$ is a manifold, we identify the unit tangent bundle
$T^1\wt M$ with $\G X$ by the map that associates to a unit tangent
vector $v\in T^1\wt M$ the unique geodesic line $\ell\in \G X$ whose
tangent vector $\dot\ell(0)$ at time $0$ is $v$. When $T^1\wt M$ is
endowed with Sasaki's metric, this mapping is an $\Isom(\wt
M)$-equivariant bi-Hölder-continuous homeomorphism, see \cite[\S
  3.1]{BroParPau19}. We denote by $\iota : \G X\ra \G X$ the
involutive {\it time-reversal map} $\ell\mapsto \{t\mapsto
\ell(-t)\}$, and again by $\iota$ the induced map from the {\it phase
  space} $\Ga\bs\G X$ to itself.  The {\it geodesic flow} on $\G X$ is
the $\Isom(X)$-equivariant one-parameter group of homeo\-morphisms
\[
\flow t :\ell\mapsto \{s \mapsto \ell(s+t)\}
\]
for all $\ell\in\G X$, with continuous time parameter $t\in\RR$ if $X$
is a manifold, and discrete time parameter $t\in\ZZ$ if $X$ is a
tree. We also denote by $(\flow t)_t$ the quotient flow on $\Ga\bs \G
X$, and call it the {\it geodesic flow} on $\Ga\bs \G X$.

A $1$-Lipschitz map $w:\RR\to X$ which is isometric on a closed
interval and has a constant value on each complementary component is a
{\it generalized geodesic}.  If $X$ is a tree, we furthermore require
that $w(0)$ and the constant values on the above-mentioned 
complementary components are vertices of $\XX$.  We denote the
Bartels-Lück metric space of generalized geodesics by $\gengeod X$,
with the metric given by Equation \eqref{eq:ggdist}, see
\cite{BarLuc12} and \cite[\S 2.2]{BroParPau19}.  The space $\gengeod
X$ isometrically contains $\G X$. The natural action of $\Isom(X)$ on
$\gengeod X$ is isometric. The geodesic flow on $\G X$ extends by the
same formula to a flow on $\gengeod X$, that we also call the {\it
  geodesic flow} on $\gengeod X$, and we use the same notation $(\flow
t)_t$ for it.

We consider a common perpendicular between two nonempty convex subsets
of $X$ of length $T$ as a generalized geodesic in $\gengeod X$, being
constant on $]-\infty,0]$ and $[T,+\infty[\,$.

\medskip
\noindent {\bf Potentials.}  We now introduce the supplementary data
(with physical origin) of potentials on $X$. Assume first that $X=\wt
M$ is a manifold. Let $\wt P:\G X=T^1\wt M\ra \RR$ be a {\it
  potential} on $X$, that is, a $\Ga$-invariant, bounded\footnote{See
\cite[\S 3.2]{BroParPau19} for a weakening of this assumption.}
Hölder-continuous real map on $T^1\wt M$.  A potential $\wt P$ is {\it
  time-reversible} if $\wt P=\wt P\circ\iota$, where $\iota$ is the
time-reversal map. For all $x,y\in \wt M$, let us define the {\it
  amplitude} $\int_x^y\wt P$ of $\wt P$ between $x$ and $y$ to be
$\int_x^y\wt P=0$ if $x=y$ and otherwise $\int_x^y\wt P=
\int_{0}^{d(x,y)} \wt P(\flow t v) \;dt$ where $v$ is the tangent
vector at $x$ to the geodesic segment from $x$ to $y$.

Now assume that $X$ is a tree, with the notation of the beginning of
Section \ref{sec:generalities}.  Denote by $V\XX$, $E\XX$ the sets of
vertices and edges of $\XX$, and by $\ov{e}$ the opposite edge of an
edge $e$.  Let $\wt c: E\XX\ra \RR$ be a (logarithmic) {\it system of
  conductances} (see for instance \cite{Zemanian91}), that is, a
$\Ga$-invariant, bounded real map on $E\XX$.  For every $\ell\in \G
X$, we denote by $e^+_0(\ell)=\ell([0,1])\in E\XX$ the first edge
followed by $\ell$, and we define the {\it potential} $\wt P:\G X\ra
\RR$ on $X$ associated with $\wt c$ as the map $\ell\mapsto \wt
c\,(e^+_0(\ell))$. The potential $\wt P$ is {\it time-reversible} if
$\wt c(\overline{e})=\wt c(e)$ for every $e\in E\XX$. For all
$x,y\in V\XX$, let $(e_1,e_2, \dots, e_k)$ be the geodesic edge path
in $\XX$ between $x$ and $y$, where $k\in\NN$ satisfies $k=0$ if and
only if $x=y$. We define the {\it amplitude} of $\wt P$ between $x$
and $y$ to be $\int_x^y\wt P= \sum_{i=1}^{k} \;\;\wt c\,(e_i)$.

If a potential $\wt P$ is time-reversible, then for all $x,y\in
\wt M$ or $x,y\in V\XX$, we have
\begin{equation}\label{eq:amplirevers}
  \int_y^x\wt P =\int_x^y\wt P\circ\iota=\int_x^y\wt P\,.
\end{equation}
We denote by $P:\Ga\bs\G X\ra\RR$ the quotient map of $\wt P$, that we
call a {\it potential} on $\Ga\bs\G X$.

Recall that every conjugacy class $[\ga]$ of loxodromic elements of
$\Ga$ defines a (not necessarily primitive) periodic orbit
$\OOO_{[\ga]}$ for the geodesic flow $(\flow t)_t$ on $\Ga\bs \G X$,
with (not necessary minimal) {\it length} $\lambda (\OOO_{[\ga]}) =
\lambda(\ga)$, so that for every $\ell\in\G X$ with $\Ga\ell\in
\OOO_{[\ga]}$, there exists a conjugate $\ga'$ of $\ga$ in $\Ga$ with
$\ell(\RR)= \Ax_{\ga'}$ and $\ell^{-1}(\ga'\ell(0))>0$. Let $\OOO$ be a
periodic orbit for the geodesic flow $(\flow t)_t$ on $\Ga\bs \G X$.
The {\it period} of $\OOO$ for the potential $\wt P$ is
\[
\int_\OOO P= \int_{\ell(0)}^{\ell(\lambda(\OOO))} \;\wt P\,,
\]
where $\ell$ is any element in $\G X$ that maps to $\OOO$ (see for
instance \cite[\S 3.1]{PauPolSha15} for the independence of the period
on the choice of $\ell$).

Let $I$ and $J$ be $\Ga$-invariant nonempty subsets of $I_\Ga$. A (not
necessarily primitive) orbit $\OOO$ of the geodesic flow on $\Ga\bs \G
X$ is called an $\{I,J\}${\em-reversible periodic orbit} if there
exists $[\ga]\in{\mathfrak R}_{I,J}$ such that $\OOO=\OOO_{[\ga]}$.
We denote by $\OOO_{I,J}$ the multiset whose underlying set
$\underline{\OOO_{I,J}}$ is the set of $\{I,J\}$-reversible periodic
orbits for the geodesic flow on $\Ga\bs \G X$, where the
multiplicities are defined by
\begin{equation}\label{eq:multperiodorbit}
  \forall\;\OOO\in \underline{\OOO_{I,J}},\quad \mult(\OOO)=
  \sum_{[\ga]\in{\mathfrak R}_{I,J}:\,\OOO=\OOO_{[\ga]}}\mult([\ga])\,.
\end{equation}
The map $\Theta_4: \underline{{\mathfrak R}_{I,J}}\ra
\underline{\OOO_{I,J}}$ defined by $[\ga]\mapsto \OOO_{[\ga]}$ is
hence a multiset bijection from ${\mathfrak R}_{I,J}$ to $\OOO_{I,J}$.

\medskip
\noindent{\bf Critical exponent and Gibbs cocycle.}  Let us fix an
arbitrary basepoint $x_*$ in $X$, which is a vertex of $\XX$ when $X$
is a tree.  The {\it critical exponent} of a potential $P$ is the
weighted (by exponential amplitudes) orbital exponential growth rate
of the group $\Ga$, defined by
\[
\delta_P= \lim_{n\ra+\infty}\;\frac{1}{n}\;\ln\;\Big(
    \sum_{\ga\in\Ga,\;n-1< d(x_*,\ga x_*)\leq n} \;\;
    \exp\big(\int_{x_*}^{\ga x_*} \wt P\;\big)\Big)\,.
\]
It is independent of the choice of the basepoint $x_*$. We have
$\delta_P\in \;]-\infty,+\infty[$ and $\delta_{P\circ\iota}=\delta_P$.

The (normalised) {\it Gibbs cocycle} of the potential $\wt P$ is the
function $C:\partial_\infty X\times \wt M\times \wt M\ra \RR$ when
$X=\wt M$ is a manifold or the function $C:\partial_\infty X\times
V\XX\times V\XX\ra \RR$ when $X$ is a tree, defined by the following
limit of difference of amplitudes for the renormalised potential
\[
(\xi,x,y)\mapsto C_\xi(x,y)= \lim_{t\ra+\infty}
\int_y^{\xi_t}(\wt P-\delta_P)-\int_{x}^{\xi_t}(\wt P -\delta_P)\,,
\]
where $t\mapsto \xi_t$ is any geodesic ray in $X$ converging to $\xi$.
The Gibbs cocycle is $\Ga$-invariant (for the diagonal action) and
locally Hölder-continuous.

\medskip
\noindent{\bf Patterson and Gibbs measures.}  A (normalised) {\it
  Patterson density} for $(\Ga,\wt P)$ is a \mbox{$\Ga$-equiv\-ariant}
family $(\mu_{x})_{x\in \wt M}$ when $X$ is a manifold, and
$(\mu_{x})_{x\in V\XX}$ when $X$ is a tree, of pairwise absolutely
continuous (positive, Borel) measures on $\partial_\infty X$, whose
support is the limit set $\Lambda\Ga$ of $\Ga$, such that
\[
\ga_*\mu_x=\mu_{\ga x}{\rm ~~~and~~~}
\frac{d\mu_x}{d\mu_y}(\xi) = e^{-C_\xi(x,\,y)}
\]
for every $\ga\in\Ga$, for all $x,y\in \wt M$ (respectively $x,y\in
V\XX$), and for (almost) every $\xi\in\partial_\infty X$.

The {\it Hopf parametrisation} of $\G X$ with basepoint $x_*$ is the
Hölder-continuous homeomorphism from $\G X$ to $(\partial_\infty X
\times \partial_\infty X\ssm{\rm Diag})\times R$, where $R=\RR$ if $X$
is a manifold and $R=\ZZ$ if $X$ is a tree, defined by $\ell\mapsto
(\ell_-,\ell_+,t)$, where $\ell_-$, $\ell_+$ are the original and
terminal points at infinity of the geodesic line $\ell$, and $t$ is
the algebraic distance along $\ell$ between the footpoint $\ell(0)$ of
$\ell$ and the closest point to $x_*$ on the geodesic line
$\ell(\RR)$.

Assume from now on that $\wt P$ is a time-reversible potential on $X$.
We denote by $dt$ the Lebesgue or counting measure on $R$. The {\it
  Gibbs measure} on $\G X$ associated with the above Patterson density
for $(\Ga,\wt P)$ is the $\sigma$-finite nonzero measure $\wt m_P$ on
$\G X$ defined using the Hopf parametrisation by
\[
d\wt m_P(\ell)=
e^{C_{\ell_-}(x_*,\,\ell(0))+C_{\ell_+}(x_*,\,\ell(0))}\;
d\mu_{x_*}(\ell_-)\;d\mu_{x_*}(\ell_+)\;dt\,.
\]
See \cite[Eq.~(4.4)]{BroParPau19}, noting that we assume $P$ to be
time-reversible. The measure $\wt m_P$ is independent of the choice of
$x_*$, is $\Ga$-invariant and $(\flow t)_{t\in R}$-invariant.
Therefore it induces\footnote{See for instance \cite[\S
  2.6]{PauPolSha15} for details on the definition of the induced
measure since $\Ga$ has torsion, hence might not act freely on $\G
X$.} a $\sigma$-finite nonzero $(\flow t)_{t\in R}\,$-invariant
measure on $\Ga\bs\G X$, called the {\it Gibbs measure} on $\Ga\bs\G
X$ for the potential $P$ and denoted by $m_P$.

\medskip
\noindent{\bf Skinning measures.}  Let $D$ be a nonempty proper closed
convex subset of $X$. We denote by $\normalpm D$ the {\it outer/inner
  unit normal bundle} of $\partial D$, that is, the space of geodesic
rays \mbox{$\rho:\pm[0,+\infty[\ra X$} with point at infinity
$\rho_\pm\in \partial_\infty X$ such that $\rho(0)\in \partial D$,
$\rho_\pm\notin \partial_\infty D$ and the closest point projection on
$D$ of $\rho_\pm$ is $\rho(0)$. We consider such a $\rho$ as a
generalised geodesic in $\gengeod X$ by requiring it to be constant on
$\mp[0,+\infty[\,$. When $X$ is a manifold and $D$ is a totally
geodesic submanifold of $X$, then the map from the unit normal bundle
$\nu^1 D$ of $D$ to $\normalout D$, sending a normal vector to $D$ to
the positive geodesic ray it defines is a homeomorphism.

Using the endpoint homeomorphism $\rho\mapsto \rho_\pm$ from
$\normalpm {D}$ to $\partial_{\infty}X \ssm\partial_{\infty}D$, we
defined in \cite{ParPau14ETDS,BroParPau19} the outer/inner {\em
  skinning measure}\footnote{This is a generalisation of the
definition of Oh and Shah \cite[\S 1.2]{OhSha12} when $D$ is a
horoball or a totally geodesic subspace in the real hyperbolic space
$\hnr$. } $\wt\sigma^\pm_{D}$ of $D$ for $(\Ga,P)$, associated with
the above Patterson density, to be the measure on $\gengeod X$, with
support contained in $\normalpm{D}$, given by
\begin{equation}\label{eq:defskin}
d\wt\sigma^\pm_D(\rho)  = 
e^{C_{\rho_\pm}(x_*,\,\rho(0))}\;d\mu_{x_*}(\rho_{\pm}) \,.
\end{equation}

Let $I'$ be an index set endowed with a left action of $\Ga$. Let
$\D=(D_i)_{i\in I'}$ be a locally finite $\Ga$-equivariant family of
nonempty proper closed convex subsets of $X$. The {\em outer/inner
  skinning measure} of $\D$ on $\gengeod X$ is the $\Ga$-invariant
locally finite measure on $\gengeod X$ defined by\footnote{See the end
of Section \ref{sec:generalities} for the definition of the
equivalence relation $\Sim[1.3]\D$. }
\[
\wt\sigma^\pm_{\D}=\sum_{i\in I'/_{\!\SimDiv\D}}\wt\sigma^\pm_{D_i}\,.
\]
It induces a locally finite measure on $\Ga \bs \gengeod X$, denoted
by $\sigma^\pm_{\D}$, called the {\em outer/inner skinning measure}
of $\D$ on $\Ga\bs \gengeod X$. Since $\wt P$ is time-reversible,
we have $\iota_*\sigma^\pm_{\D}=\sigma^\mp_{\D}$, and in particular
the measures $\sigma^-_{\D}$ and $\sigma^+_{\D}$ have the same total
mass:
\begin{equation}\label{eq:sameskinmass}
\|\sigma^-_{\D}\|=\|\sigma^+_{\D}\|\,.
\end{equation}

\section{Counting strongly reversible closed geodesics}
\label{sec:proofs}

In this Section, let $(X,\Ga)$ be as in the beginning of Section
\ref{sec:generalities}.  Let $x_*$ be an arbitrary basepoint in $X$,
that we take to be a vertex of $\XX$ when $X$ is a tree.

Let $I$ and $J$ be two nonempty $\Ga$-invariant subsets of $I_\Ga$.
Let $\F_I= (F_\alpha)_{\alpha\in I}$ (respectively $\F_J= (F_\alpha)
_{\alpha\in J}$) be the $\Ga$-equivariant family of the fixed point
sets of the elements of $I$ (respectively $J$). Let $\sigma^+_{I} =
\sigma^+_{\F_I}$ be the outer skinning measure of $\F_I$, and
$\sigma^-_{J} = \sigma^-_{\F_J}$ the inner skinning measure of $\F_J$.
Let $\wt P:\G X\ra \RR$ be a time-reversible potential with positive
critical exponent $\delta_P$. When $X$ is a tree, we furthermore
assume that the smallest nonempty $\Ga$-invariant simplicial subtree
of $\XX$ is uniform,\footnote{This means that there exists a discrete
subgroup of $\Isom(X)$ with compact quotient.}  without vertices
of degree $2$.

When $X$ is a tree, we denote by $L_\Ga$ the {\it length spectrum} of
$\Ga$, that is, the subgroup of $\ZZ$ generated by the translation
lengths in $X$ of the loxodromic elements of $\Ga$. By
\cite[Lem.~4.18]{BroParPau19}, under the assumption of the finiteness
of the measure $m_P$, the above assumption on the tree $\XX$ implies
that $L_\Ga$ is either $\ZZ$ or $2\ZZ$. When $L_\Ga=\ZZ$, the geodesic
flow $(\flow{t})_{t\in\ZZ}$ is mixing for the Gibbs measure $m_{P}$ on
$\Ga\bs \G X$ by \cite[Theo.~4.17]{BroParPau19}. Let $\gengeod_{\rm
  even} X$ be the subset of elements $\ell\in \gengeod X$ such that
${\displaystyle \lim_{t\ra-\infty}}\,\ell(t)$, $\ell(0)$,
${\displaystyle \lim_{t\ra+\infty}}\,\ell(t)$ are at even distance of
$x_*$ or belong to $\partial_\infty X$. Let
\[
\G_{\rm even} X=\G X\cap \gengeod_{\rm even} X
=\{\ell\in\G X: d(x_*,\ell(0)) \in 2\ZZ\}\,.
\]
By loc.~cit., when $L_\Ga =2\ZZ$, the subsets $\gengeod_{\rm even} X$
and $\G_{\rm even} X$ are $\Ga$-invariant, and the even times geodesic
flow $(\flow{2t})_{t\in\ZZ}$ is mixing for the restriction $m_{P,\,
  {\rm even}}$ to $\Ga\bs \G_{\rm even} X$ of the Gibbs measure
$m_{P}$. We also denote by $\sigma^+_{I,\,{\rm even}}$ and
$\sigma^-_{J,\,{\rm even}}$ the restriction of $\sigma^+_I$ and
$\sigma^+_J$ to $\Ga\bs \gengeod_{\rm even} X$.

Let $\delta_P^*= \delta_P$ when $X$ is a manifold. When $X$ is a tree,
let $\delta_P^*=1- e^{-\delta_P}$ when $L_\Ga=\ZZ$ and $\delta_P^*=
1-e^{-2\delta_P}$ when $L_\Ga=2\ZZ$. For every $s>0$ when $X$ is a
manifold, and every $s\in\NN\ssm\{0\}$ when $X$ is a tree with
$L_\Ga=\ZZ$, let $\Perp^*(\F_I,\F_J,s) = \Perp(\F_I,\F_J,s)$,
$IJ^*(s)=IJ(s)$ and ${\mathfrak R} ^*_{I,J} (s)={\mathfrak R}_{I,J}
(s)$.

When $X$ is a tree with $L_\Ga=2\ZZ$, for every $s\in 2\NN\ssm\{0\}$,
let $\Perp^*(\F_I,\F_J,s)$ be the multiset restriction of the common
perpendiculars in $\Perp(\F_I,\F_J,s)$ with both endpoints (which are
vertices of $\XX$) at even distance from the basepoint $x_*$. Note
that when the fixed point sets of the elements in $I\cup J$ are at
even distance from the base point $x_*$, by the relations between the
three distances between three points in a tree, we have
$\Perp^*(\F_I,\F_J,s) = \Perp(\F_I, \F_J,s)$. Let also $IJ^*(s)$ be
the multiset restriction of $IJ(s)$ consisting of its elements
$\Ga(\alpha,\beta)$ with both endpoints of the common perpendicular
$\wt\rho_{\alpha,\beta}$ at even distance from $x_*$.  Finally let
${\mathfrak R}^*_{I,J}(2s)$ be the image in ${\mathfrak R}_{I,J} (2s)$
of $IJ^*(s)$ by the multiset bijection $\Theta_2$ defined in
Proposition \ref{prop:bijcomperprecloxo}. Note that the multiset
bijection $\Theta_1$ defined in that proposition sends $IJ^*(s)$ to
$\Perp^*(\F_I,\F_J,s)$.

Recall the following standard convention on sums over multisets. Let
$S=(\underline S,\mult)$ be a finite real-valued multiset.  For every
map $f$ from $\underline{S}$ to a real vector space, let
$$\sum_{a\in S}
f(a)=\sum_{a\in\underline{S}}\; \mult(a)\,f(a)\,.
$$

Let 
\begin{equation}\label{eq:defiNIPT}
  \N_{I,J,P}(T)=\sum_{[\ga]\,\in \,{\mathfrak R}^*_{I,J}(T)}
  \;e^{\frac{1}{2}\int_{[\ga]}P}
\end{equation}
be the counting function of conjugacy classes of $\{I,J\}$-reversible
loxodromic elements of $\Ga$ with translation length at most $T$, with
multiplicities given by Equation \eqref{eq:multreclox} and with
weights given by their exponential half-periods

\btheo\label{theo:main} Let $(X,\Ga,\wt P\,)$ be as above.

\hypertarget{theomaincount1}{(1)} Assume that $X$
is a manifold or a tree with $L_\Ga=\ZZ$. Assume that the Gibbs
measure $m_{P}$ is finite and mixing under the geodesic flow on
$\Ga\bs \G X$, and that the skinning measures $\sigma^+_{I}$ and
$\sigma^-_{J}$ are nonzero and finite. As $T\ra+\infty$ with $T\in
2\ZZ$ if $X$ is a tree, we have
\begin{equation}\label{eq:equivmain}
\N_{I,J,P}(T)\sim
\frac{\|\sigma^+_{I}\|\;\|\sigma^-_{J}\|}{\delta_P^*\;\|m_P\|}\,
\exp\big(\,\frac{\delta_P}{2}\;T\,\big)\,. 
\end{equation}
Furthermore, if $\wt P=0$ and either $X=\wt M$ is a symmetric space,
$\Ga\bs \wt M$ has finite volume and exponentially mixing geodesic
flow, or if $X$ is a tree and $\Ga$ is a geometrically finite tree
lattice with $L_\Ga =\ZZ$, then there is an additive error term of the
form $\bigO(e^{(\frac{\delta_P}{2} -\kappa)T})$ for some $\kappa>0$ in
Equation \eqref{eq:equivmain}.

\hypertarget{theomaincount2}{(2)} Assume that $X$
is a tree with $m_P$ finite and $L_\Ga=2\ZZ$, and that the skinning
measures $\sigma^+_{I,\,{\rm even}}$ and $\sigma^-_{J,\,{\rm even}}$
are nonzero and finite.  As $T\ra+\infty$ with $T\in 4\ZZ$, we have
\begin{equation}\label{eq:equivmain2Z}
  \N_{I,J,P}(T)\sim \frac{\|\sigma^+_{I,\,{\rm even}}\|\;
\|\sigma^-_{J,\,{\rm even}}\|}{\delta_P^*\;\|m_{P,\,{\rm even}}\|}\,
\exp\big(\,\frac{\delta_P}{2}\;T\,\big)\,.
\end{equation}
Furthermore, if $\wt P=0$ and $\Ga$ is a geometrically finite tree
lattice with $L_\Ga=2\ZZ$, then there is an additive error term of the
form $\bigO(e^{(\frac{\delta_P}{2} -\kappa)T})$ for some $\kappa>0$ in
Equation \eqref{eq:equivmain2Z}.  \etheo

When $X=\wt M$ is a manifold, the mixing assumption for the error term
in Theorem \ref{theo:main} (\hyperlink{theomaincount1}{1}) holds for
instance if $\wt M$ is a real hyperbolic space $\HH^n_\RR$ and $\Ga$
is geometrically finite by \cite{LiPan22}, or if $\Ga$ is an
arithmetic lattice in the isometry group of a negatively curved
symmetric space $\wt M$ by the works of Clozel, Kleinbock and
Margulis, as explained in \cite[\S 9.1]{BroParPau19}. When $X$ is a
tree, the error term in Theorem \ref{theo:main}
(\hyperlink{theomaincount1}{1}) (respectively
(\hyperlink{theomaincount2}{2})) holds if $X$ is the Bruhat-Tits tree
of a rank one simple algebraic group $G$ over a nonarchimedean local
field and $\Ga$ is any lattice in $G$ with $L_\Ga=\ZZ$ (respectively
with $L_\Ga=2\ZZ$) by \cite[Theo.~C]{Lubotzky91}.

Theorem \ref{theo:mainintro} (\hyperlink{mainintro1}{1}) in the
introduction follows from Theorem \ref{theo:main}
(\hyperlink{theomaincount1}{1}) by restricting to the manifold case
and by taking $\wt P=0$ (which is time-reversible), in which case
$\delta_P=\delta_\Ga>0$ and $m_P=m_{\rm BM}$, since $\N_{I,J}(T)$ as
defined in the introduction is equal to $\N_{I,J,0}(T)$ as defined in
Equation \eqref{eq:defiNIPT}.

\medskip
\dem The main idea of the proof is to relate the counting function of
strongly reversible closed geodesics in $\Ga\bs X$ with the counting
function of the common perpendiculars between the pairs of fixed point
sets of the involutions in $I$ and $J$, using the preliminary work of
Section \ref{sec:generalities}.

Let $(\alpha,\beta) \in \wt{IJ}$. Let $\ga=\beta\alpha$ and let
$\rho_{\alpha,\beta}$ be the $\Ga$-orbit of the common perpendicular
$[x,y]$ from $F_\alpha$ to $F_\beta$ with $x\in F_\alpha$, so that
$\beta y=y$. Let us define $\int_{\rho_{\alpha,\beta}} P=\int_x^y\wt
P$, which is well defined since $\wt P$ is $\Ga$-invariant. Since $\wt
P$ is time-reversible, we have $\int_x^y\wt P= \int_y^{x}\wt
P\circ\iota= \int_y^{x}\wt P$. Hence
\begin{equation}\label{eq:reversamplitude}
\int_{\rho_{\alpha,\beta}} P=\int_{\rho_{\beta,\alpha}} P\,.
\end{equation}
Since $\wt P$ is $\Ga$-invariant, we have $\int_x^y\wt P=\int_y^x \wt
P=\int_{\beta y}^{\beta x}\wt P= \int_y^{\beta x}\wt P$. Hence since
$x$ belongs to the translation axis of $\ga$ by Lemma
\ref{lem:prodinversion}, since $y$ is the midpoint of $[x,\ga x]=
[x,\beta x]$, and by the additivity properties of the amplitudes, we
have
\begin{equation}\label{eq:relatperiod}
  \int_{[\ga]} P=\int_x^{\ga x}\wt P=\int_x^{y}\wt
  P+\int_y^{\beta x}\wt P =2\int_{\rho_{\alpha,\beta}} P\,.
\end{equation}

The assumptions of \cite[Coro.~12.3]{BroParPau19} when $X$ is a
manifold, and of \cite[Theo.~12.9]{BroParPau19} when $X$ is a tree
with $L_\Ga=\ZZ$, are satisfied by the assumptions of Theorem
\ref{theo:main} (\hyperlink{theomaincount1}{1}). The assumptions of
\cite[Theo.~12.12]{BroParPau19} when $X$ is a tree with $L_\Ga=2\ZZ$
are satisfied by the assumptions of Theorem \ref{theo:main}
(\hyperlink{theomaincount2}{2}) and by the second part of
\cite[Theo.~4.17]{BroParPau19} for the required mixing property. We
claim that the counting parts in these results, with the above
definitions of $\delta_P^*$ and $\Perp^*$, imply that
\begin{equation}\label{eq:countcomperp}
  \sum_{\rho\,\in\, \Perp^*(\F_I,\F_J,\frac{T}{2})} \;e^{\int_{\rho}P}\sim
\frac{\|\sigma^+_{I}\|\;\|\sigma^-_{J}\|}{\delta_P^*\,\|m_P\|}\;
\exp\big(\,\delta_P\;\frac{T}{2}\,\big)\,,
\end{equation}
as $T\ra+\infty$ with $T\in \RR$ if $X$ is a manifold and $T\in 2\ZZ$
if $X$ is a tree with $L_\Ga=\ZZ$, and that
\begin{equation}\label{eq:countcomperpeven}
\sum_{\rho\,\in\, \Perp^*(\F_I,\F_J,\frac{T}{2})}\;e^{\int_{\rho}P}\sim
  \frac{\|\sigma^+_{I,\,{\rm even}}\;\|\sigma^-_{J,\,{\rm even}}\|}
       {\delta_P^*\;\|m_{P,\,{\rm even}}\|}\;
       \exp\big(\,\delta_P\;\frac{T}{2}\,\big)\,,
\end{equation}
as $T\ra+\infty$ with $T\in 4\ZZ$ if $X$ is a tree with $L_\Ga=2\ZZ$.

The conclusion of \cite[Coro.~12.3]{BroParPau19} is formulated in a
way that is not directly applicable, considering the difference
between the counting function $\N_{\F_I,\F_J,P}$ in \cite[page 255]
{BroParPau19} and the above counting function $\sum_{\rho\,\in\,
  \Perp(\F_I,\F_J, \frac{T}{2})} \;e^{\int_{\rho}P}$.  The conclusions
of \cite[Theo.~12.9 and Theo.~12.12]{BroParPau19} would only apply
when $\Ga\bs I$ and $\Ga\bs J$ are reduced to one element (that is,
when $I$ and $J$ contain only one conjugacy class of involutions, as
for Corollary \ref{cor:introNagao}). Hence we need to go over the
proofs of these results and explain the work on the multiplicities in
order to obtain Equations \eqref{eq:countcomperp} and
\eqref{eq:countcomperpeven}, without stating in full these statements
and their proofs.

The process is the following one. The counting part in these three
statements \cite[Coro.~12.3, Theo.~12.9 and Theo.~12.12]{BroParPau19}
are deduced from three claims of narrow convergence of measures in the
quotient space $(\Ga\times\Ga) \bs(\gengeod X\times \gengeod X)$ by
integrating on the constant function $1$. These three claims are
themselves consequences of three convergence statements of measures on
the product space $\gengeod X\times \gengeod X$ stated as
\cite[Theo.~11.1, Theo.~11.9 and Eq.~(11.28)]{BroParPau19}.  These
last three convergence statements have the form
\begin{equation}\label{eq:weakstarup}
\lim_{t\ra+\infty}\sum_{a\in \A_t} \;\wt\nu_a\;=\;\wt\nu_\infty\,,
\end{equation}
where

$\bullet$~ for every $t>0$, the index set $\A_t$ is endowed with a left
action of $\Ga\times\Ga$,

$\bullet$~ the measures $\wt\nu_a$ for $a\in \A_t$ are positive
multiples of Dirac masses on $\gengeod X\times \gengeod X$ and the map
$a\mapsto \wt\nu_a$ is equivariant for the action of $\Ga\times\Ga$.

In order to be continuous, the process of defining measures on the
quotient by branched coverings (see for instance \cite[\S 2.6]
{PauPolSha15}) consists of the following operations. The new index set
of the sum is the quotient modulo the action of the group $\Ga\times
\Ga$ on the index set $\A_t$.  The image measure $\nu_a$ of $\wt\nu_a$
in $(\Ga\times\Ga)\bs(\gengeod X\times\gengeod X)$ is divided by the
cardinality of the stabiliser in $\Ga\times\Ga$ of the support of
$\wt\nu_a$. Let $\nu_\infty$ be the measure on $(\Ga\times\Ga)\bs
(\gengeod X\times\gengeod X)$ induced by $\wt\nu_\infty$. This gives a
convergence, which turns out to be a narrow convergence as explained
in \cite[Coro.~12.3, Theo.~12.9 and Theo.~12.12]{BroParPau19}, of the
form
\begin{equation}\label{eq:weakstardown}
\lim_{t\ra+\infty}\sum_{a\in (\Ga\times\Ga)\bs \A_t}\;
\frac{1}{\operatorname{Card} (\operatorname{Stab}_{\Ga\times\Ga}(a))}
\;\nu_a =\nu_\infty\,.
\end{equation}

Let us now explain how we apply this process. For every $\alpha\in I$
(resp.~$\beta\in J$), let $\ov\alpha$ (resp.~$\ov\beta$) be its class
in $I/_{\sim_I}$ (resp.~ $J/_{\sim_J}$). When $X$ is a manifold, the
statement in \cite[Theo.~11.1]{BroParPau19} (taking into account that
the action of $\Ga$ on $I$ is by conjugation) uses the index set
\[
\A_t=\big\{(\ov\alpha,\ov\beta,\ga)\in
I/_{\sim_I}\times J/_{\sim_J}\times\Ga:
0<d(F_\alpha,F_{\ga\beta\ga^{-1}})\leq t\big\}\,
\]
where $t$ varies in $]0,+\infty[\,$. The action of $\Ga\times\Ga$ on
$\A_t$ is given by
\[
(\ga_1,\ga_2)\cdot
(\ov\alpha,\ov\beta,\ga) =(\ga_1\,\ov\alpha\,\ga_1^{-1}, \ga_2\,
\ov\beta\, \ga_2^{-1},\ga_1\,\ga\,\ga_2^{-1})\,.
\]
If $(\ga_1,\ga_2)\in \Ga\times\Ga$ fixes $(\ov\alpha,\ov\beta,\ga)\in
\A_t$, then $\ga_2= \ga^{-1} \ga_1\ga$ is uniquely determined by
$\ga_1$. Note that $\ga_2=\ga^{-1} \ga_1\ga$ belongs to
$\Ga_{F_\beta}$ if and only if $\ga_1$ belongs to $\ga \Ga_{F_\beta}
\ga^{-1}= \Ga_{\ga F_\beta} =\Ga_{F_{\ga\beta\ga^{-1}}}$. For every
$(\alpha_0,\ga_0)\in I\times\Ga$, by the definition of the equivalence
relation $\sim_I$ in Formula \eqref{def:simi},
we have $\alpha_0 \sim_I \ga_0\,\alpha_0\, \ga_0^{-1}$ if and only if
$F_{\alpha_0} = F_{\ga_0\,\alpha_0\,\ga_0^{-1}}$, that is, if and only
if $\ga_0$ belongs to the stabilizer $\Ga_{F_{\alpha_0}}$ of
$F_{\alpha_0}$ in $\Ga$. Applying this for both $(\alpha_0,\ga_0)
=(\alpha,\ga_1)$ and $(\alpha_0,\ga_0)=(\beta,\ga_2)$ with $I$
replaced by $J$, this proves that the stabilizer of $(\ov\alpha,
\ov\beta,\ga)\in \A_t$ in $\Ga \times\Ga$ is isomorphic to
$\Ga_{F_\alpha}\cap\Ga_{F_{\ga\beta\ga^{-1}}}$, by the map $(\ga_1,
\ga_2) \mapsto \ga_1$. Hence, these two finite groups have the same
cardinality.

For every $t\in[0,+\infty[\,$, let us consider the set
\[
\A'_t=\big\{(\ov\alpha,\ov\beta)\in I/_{\sim_I} \times
J/_{\sim_J}: 0<d(F_\alpha,F_{\beta})\leq t\big\}
\]
endowed with the diagonal action $\ga\cdot(\ov\alpha,\ov\beta)=
(\ga\,\ov\alpha\,\ga^{-1},\ga\,\ov\beta\,\ga^{-1})$ by conjugation of
$\Ga$. The map $\wt \Theta: (\ov\alpha,\ov\beta)\mapsto(\ov\alpha,
\ov\beta,e)$ from the set $\A'_t$ endowed with its action of $\Ga$ to
the set $\A_t$ endowed with its action of $\Ga\times\Ga$ is
equivariant for the diagonal group morphism from $\Ga$ into
$\Ga\times\Ga$. It hence induces, by taking quotients, a map
\[
\Theta:\Ga \bs \A'_t\ra(\Ga\times\Ga)\bs \A_t\,.
\]
Note that for all $\alpha\in I$, $\beta\in J$ and $\ga\in\Ga$, we have
$(e,\ga) \cdot (\ov\alpha, \ov\beta,\ga)= (\ov\alpha, \ga\,\ov\beta\,
\ga^{-1},e)$. Hence, every element of $\A_t$ is in the same $(\Ga
\times \Ga)$-orbit as an element of the image of $\wt\Theta$, and the
map $\Theta$ is surjective. The map $\Theta$ is also injective, since
for all $(\ov\alpha,\ov\beta), (\ov{\alpha'}, \ov{\beta'})\in \A'_t$,
if there exists $(\ga_1,\ga_2)\in\Ga\times\Ga$ such that $(\ga_1,
\ga_2) \cdot (\ov\alpha,\ov\beta,e) = (\ov{\alpha'}, \ov{\beta'},e)$,
then $\ga_1=\ga_2$, and $(\ov\alpha,\ov\beta)$ and $(\ov{\alpha'},
\ov{\beta'})$ are in the same $\Ga$-orbit in $\A'_t$.

The map $\Theta':\Ga (\alpha,\beta)\mapsto \Ga(\ov\alpha,\ov\beta)$
from the set $\underline{IJ(t)}$ (defined in Equation
\eqref{eq:defiIIGaT}) to the set $\Ga \bs \A'_t$ is surjective. Its
fiber over $\Ga(\ov\alpha, \ov\beta)$ has cardinality
$|\Ga_{F_\alpha}/Z_\Ga(\alpha)|\; |\Ga_{F_\beta}/ Z_\Ga(\beta)|$ by
Lemma \ref{lem:cardequivmodI}.  Hence

$\bullet$~ by the definitions of a sum over a multiset and of the
multiplicity of a common perpendicular between elements of $\F_I$ and
$\F_J$ given in Equation \eqref{eq:multcomperp} for the first equality
below,

$\bullet$~ by the orders of the fibers of $\Theta'$ for the
second equality below,

$\bullet$~ using the bijection $\Theta$ for the
third equality below,

\noindent we have
\begin{align*}
  \sum_{\rho\,\in\, \Perp(\F_I,\,\F_J,\,t)} \;e^{\int_{\rho}P}&=
  \sum_{\Ga(\alpha,\beta)\,\in\, \underline{IJ(t)}}
  \frac{e^{\int_{\rho_{\alpha,\beta}}P}}{\card(\Ga_{F_\alpha}\cap\Ga_{F_\beta})\;
    |\Ga_{F_\alpha}/Z_\Ga(\alpha)|\; |\Ga_{F_\beta}/Z_\Ga(\beta)|}
  \\&=\sum_{\Ga\cdot(\ov\alpha,\ov\beta)\,\in\, \Ga \bs \A'_t}
  \frac{e^{\int_{\rho_{\alpha,\beta}}P}}{\card(\Ga_{F_\alpha}\cap\Ga_{F_\beta})}
\\&=\sum_{(\Ga\times\Ga)\cdot(\ov\alpha,\ov\beta,\ga)\,\in\, (\Ga\times\Ga) \bs \A_t}
\frac{e^{\int_{\rho_{\alpha,\ga\beta\ga^{-1}}}P}}
     {\card(\Ga_{F_\alpha}\cap\Ga_{F_{\ga\beta\ga^{-1}}})}\,.
\end{align*}
Therefore, taking $t=\frac{T}{2}$, Equation \eqref{eq:countcomperp}
when $X$ is a manifold (so that $\delta_P^*= \delta_P$) then follows
from \cite[Theo.~11.1]{BroParPau19} applied with $\D^-=\F_I$ and
$\D^+=\F_J$, by the process explained above starting from Equation
\eqref{eq:weakstarup}, passing through Equation
\eqref{eq:weakstardown} and integrating on the constant function $1$.

When $X$ is a tree and $L_\Ga=\ZZ$, the proof of Equation
\eqref{eq:countcomperp} is similar, taking $t=\frac{T}{2}\in\NN$
(hence $T\in 2\NN$), replacing \cite[Theo.~11.1]{BroParPau19} by
\cite[Theo.~11.9]{BroParPau19} applied with $\D^-=\F_I$ and
$\D^+=\F_J$, using the above definition of
$\delta_P^*=1-e^{-\delta_P}$ and the definition of the amplitudes for
trees given in Section \ref{sec:dyn}.

When $X$ is a tree and $L_\Ga=2\ZZ$, the proof of Equation
\eqref{eq:countcomperpeven} is also similar, taking
$t=\frac{T}{2}\in2\NN$ (hence $T\in 4\NN$), replacing
\cite[Theo.~11.1]{BroParPau19} by \cite[Eq.~(11.28)]{BroParPau19}
applied with $\D^-=\F_I$ and $\D^+=\F_J$, using the above definitions
of $\delta_P^*=1-e^{-2\delta_P}$ and $\Perp^*$, noting that we have
$m_{P,\, {\rm even}}= \frac{\|m_P\|}{2}$.

\medskip
Let us prove Assertion (\hyperlink{theomaincount1}{1}) of Theorem
\ref{theo:main}. Under its assumptions, using respectively in the
following sequence of equalities and equivalences

$\bullet$~ the definition of $\N_{I,J,P}(T)$ in Equation
\eqref{eq:defiNIPT} for the first equality,

$\bullet$~ the fact that the map $\Theta_2$ in Proposition
\ref{prop:bijcomperprecloxo} is a multiset bijection sending
$IJ(\frac{T}{2})$ to ${\mathfrak R}_{I,J}(T)$, and Equation
\eqref{eq:relatperiod}, for the second equality,

$\bullet$~ the fact that the map $\Theta_1$ in Proposition
\ref{prop:bijcomperprecloxo} is a multiset bijection sending
$IJ(\frac{T}{2})$ to $\Perp(\F_I,\F_J,\frac{T}{2})$, for the third
equality, and

$\bullet$~ Equation \eqref{eq:countcomperp} for the final equivalence,

\noindent as $T\ra+\infty$ with $T\in 2\ZZ$ if $X$ is a tree with
$L_\Ga=\ZZ$, we
have
\begin{align}
  \N_{I,J,P}(T)&=\sum_{[\ga]\,\in \,{\mathfrak R}^*_{I,J}(T)}
  e^{\frac{1}{2}\int_{[\ga]}P}=
  \sum_{\Ga(\alpha,\beta)\,\in \,IJ^*(\frac{T}{2})}
  e^{\int_{\rho_{\alpha,\beta}}P}\nonumber\\&=
  \sum_{\rho\,\in\, \Perp^*(\F_I,\F_J,\frac{T}{2})}
  e^{\int_{\rho}P}\sim
\frac{\|\sigma^+_{I}\|\,\|\sigma^-_{J}\|}{\delta_P^*\,\|m_P\|}\,
\exp\big(\,\frac{\delta_P}{2}\;T\,\big)\,.\label{eq:relatNINPerp}
\end{align}
This proves Equation \eqref{eq:equivmain}.

When $X=\wt M$ is a manifold, the final error term claim of Theorem
\ref{theo:main} (\hyperlink{theomaincount1}{1}) follows from the
analogous error term in Equation \eqref{eq:countcomperp}, obtained in
\cite[Theo.~12.7 (2)]{BroParPau19} using

$\bullet$~ the fact that when $\wt P=0$, $\wt M$ is locally symmetric
and $\Ga\bs\wt M$ has finite volume, then $m_P=m_{\rm BM}$ is the
(smooth) Liouville measure, which is finite,

$\bullet$~ the fact that the fixed point sets of involutions of $\Ga$
are (smooth) totally geodesic submanifolds of $\wt M$,

$\bullet$~ the exponential mixing assumption that is needed in
\cite[Theo.~12.7 (2)]{BroParPau19}.

When $X$ is a tree with $L_\Ga=\ZZ$, the final error term claim of
Theorem \ref{theo:main} (\hyperlink{theomaincount1}{1}) follows from
the analogous error term in Equation \eqref{eq:countcomperp}, now
obtained by \cite[Theo.~12.16]{BroParPau19}, and more precisely in its
following Remark (i) page 281.

\medskip The proof of Assertion (\hyperlink{theomaincount2}{2}) of
Theorem \ref{theo:main} when $X$ is a tree with $L_\Ga=2\ZZ$ is
similar. Note that by the above definitions of $IJ^*(T)$,
$\Perp^*(\F_I,\F_J,T)$ and ${\mathfrak R}^*_{I,J}(T)$, and as
mentioned before Equation \eqref{eq:defiNIPT}, the map $\Theta_1$
sends $IJ^*(T)$ to $\Perp^*(\F_I,\F_J, T)$ and $\Theta_2$ sends
$IJ^*(T)$ to ${\mathfrak R}^*_{I,J}(2T)$. We then obtain Equation
\eqref{eq:relatNINPerp} when $T\in 4\ZZ$ by replacing Equation
\eqref{eq:countcomperp} by Equation \eqref{eq:countcomperpeven}. For
the error term, we replace Remark (i) following
\cite[Theo.~12.16]{BroParPau19} by its following Remark (ii).
\cqfd

\section{Equidistribution of strongly reversible closed geodesics}
\label{sec:eqproof}

Let $(X,\Ga)$ be as in the beginning of Section
\ref{sec:generalities}. Let $I$ and $J$ be $\Ga$-invariant nonempty
subsets of $I_\Ga$. Let $\wt P$ be a time-reversible potential on $X$
with positive critical exponent $\delta_P$.

In this Section, we prove that the Lebesgue measures along the
$\{I,I\}$-reversible periodic orbits of the geodesic flow on $\Ga\bs\G
X$ of length at most $T$, with multiplicities and with weights given
by their exponential half-periods, equidistribute towards the Gibbs
measure $m_P$, see Theorem \ref{theo:main2}.  This result is due to
\cite[Thm.\,1.4]{ErlSou24} in the case of surfaces with constant
negative curvature and orientation-preserving involutions with $\wt
P=0$. Here, we improve even this result by adding an error term for
the equidistribution.

Note also that the collection of periodic orbits considered in Theorem
\ref{theo:main2} is considerably smaller than the full collection of
periodic orbits that is known to equidistribute to the Gibbs measure
by \cite[Thm.\, 9.11]{PauPolSha15}, generalizing results of
\cite{Bowen72} and \cite{Roblin03}.

We only consider the case when $X$ is a manifold for simplicity. The
analogous result when $X$ is a tree should be true, but this would
require a complete tree version of \cite{ParPau25b}. This has only be
done in \cite[\S 6]{ParPauSay25c} for common perpendiculars between
horoballs and when $\wt P=0$.

Let $\ga\in \wt{\mathfrak R}_{I,J}$.  We fix a point $x_\ga\in
\Ax_\ga$. We parametrize the translation axis of $\ga$ by the unique
isometric map $\Ax_\ga: \RR\ra X$ such that $\Ax_\ga(0)=x_\ga$ and
$\Ax_\ga(\lambda(\ga)) = \ga\,x_\ga$.  We denote by $\Leb_{[\ga]}$ the
measure on $\Ga\bs\G X$ obtained by pushing forward the Lebesgue
measure on the real interval $[0,\lambda(\ga)]$ by the map $t\mapsto
\Ga\flow{t}\Ax_\ga$. This measure does not depend on the choice of
$x_\ga$ by the invariance of the Lebesgue measure by translations, nor
on the choice of representative $\ga$ of the conjugacy class $[\ga]$.

For every periodic orbit $\OOO$ of the geodesic flow on $\Ga \bs\G X$,
with $[\ga]$ any conjugacy class in $\Ga$ such that $\OOO
=\OOO_{[\ga]}$, we define $\Leb_\OOO=\Leb_{[\ga]}$, which does not
depend on the above choice of $[\ga]$.

For every $(\alpha,\beta)\in\wt{IJ}$, let $\wh{\rho}\,_{\alpha,\beta}
\in\G X$ be the unique geodesic line such that $\wh{\rho}
\,_{\alpha,\beta}$ coincides with $\wt{\rho}\,_{\alpha ,\beta}$ on
$[0,\lambda(\rho) ]$. We denote by $\Leb_\rho$ the measure on
$\Ga\bs\G X$ obtained by pushing forward the Lebesgue measure on the
real interval $[0, \lambda(\rho)]$ by the map $t\mapsto \Ga\flow{t}
\wh{\rho}\,_{\alpha ,\beta}$, where Note that in general, the two
measures $\Leb_{\rho_{\alpha,\beta}}$ and $\Leb_{\rho_{\beta,\alpha}}$
do not coincide: Otherwise, by taking the footpoint map on their
support, we would have $\rho_{\alpha,\beta} =\rho_{\beta,\alpha}$.
There would hence exist $\ga\in\Ga$ such that $\ga\,\wt
\rho_{\beta,\alpha} =\wt\rho_{\alpha,\beta}$. But such an element
$\ga$ would fix the midpoint of the geodesic segment
$\wt\rho_{\beta,\alpha}$, and there is no reason why such an elliptic
element would always exist. This is the reason why we will only
consider $I=J$ in Theorem \ref{theo:main2}.

\blemm If $X$ is a manifold, for every element $\Ga (\alpha, \beta)\in
IJ$, we have
\begin{equation}\label{eq:decompLeb}
\Leb_{\Theta_2(\Ga(\alpha,\beta))} =
\Leb_{\Theta_1(\Ga(\alpha,\beta))}+\Leb_{\Theta_1(\Ga(\beta,\alpha))}\,.
\end{equation}
\elemm

\dem Let $(\alpha,\beta)\in \wt{IJ}$ and $\ga=\beta\alpha$, so that
with $\Theta_1$ and $\Theta_2$ the multiset bijections defined in
Proposition \ref{prop:bijcomperprecloxo}, we have $\Theta_1
(\Ga(\alpha, \beta))= \rho_{\alpha,\beta}$ and $\Theta_2
(\Ga(\alpha,\beta))=[\ga]$.  Let $x_\alpha$ and $x_\beta$ be the
endpoints of the common perpendicular $\wt\rho_{\alpha,\beta}$ with
$x_\alpha\in F_\alpha$ and $x_\beta\in F_\beta$. Note that the
endpoints of the common perpendicular $\wt\rho_{\beta,\beta\alpha
  \beta^{-1}}$ are $x_\beta\in F_\beta$ and $\beta x_\alpha\in \beta
F_\alpha=F_{\beta\alpha\beta^{-1}}$.  Let $s= d(x_\alpha, x_\beta)$,
so that the translation length of $\ga$ is $2s$ by Equation
\eqref{eq:length2length}. Let $\tau_s: t\mapsto t+ s$ be the
translation by $s$ on the real line $\RR$. Let $\Ax_\ga:\RR \ra X$
with $\Ax_\ga(0)=x_\ga$ and $\wh\rho_{\alpha, \beta} :\RR\ra X$ be the
geodesic lines defined above.  By construction and Equation
\eqref{eq:equivarthoalphabeta}, we have
\[
\Ax_\ga=\wh\rho_{\alpha,\beta}\quad \text{and}\quad
\flow{s}\,\wh\rho_{\alpha,\beta}=\wh\rho_{\beta,\beta\alpha\beta^{-1}}=
\beta\,\wh\rho_{\beta,\alpha}\,.
\]
Hence the maps $\Psi:\RR\ra \Ga\bs\G X$ defined by $t\mapsto \Ga
\flow{t} \Ax_\ga=\Ga\flow{t}\wh\rho_{\alpha,\beta}$ and $\Psi':\RR\ra
\Ga \bs\G X$ defined by $t\mapsto \Ga\flow{t}\wh\rho_{\beta,\alpha}$
satisfy $\Psi\circ \tau_s=\Psi'$.  Let $\Leb$ be the Lebesgue measure
on $\RR$. Using in the following sequence of equalities

$\bullet$~ the definition of $\Leb_{[\ga]}$ for the first equality,

$\bullet$~ the fact that $\Leb\restrict{ [0,2s]}=\Leb\restrict{ [0,s]}
+\Leb\restrict{ [s,2s]}$ and the linearity of $\Psi_*$ for the second
equality,

$\bullet$~ the fact that $\Psi\circ \tau_s=\Psi'$ for the third equality,

$\bullet$~ the definition of $\Leb_{\rho_{\alpha,\beta}}$ for the last equality,

\noindent we have
\begin{align}
  \Leb_{[\ga]}&=\Psi_*(\Leb\restrict{
    [0,2s]})= \Psi_*(\Leb\restrict{ [0,s]})+\Psi_*(\Leb\restrict{
    [s,2s]})\nonumber\\& =\Psi_*(\Leb\restrict{ [0,s]})+
  \Psi'_*(\Leb\restrict{ [0,s]})= \Leb_{\rho_{\alpha,\beta}}+
  \Leb_{\rho_{\beta,\alpha}}\,.\label{eq:detailsdecompLeb}
\end{align}
This proves Equation \eqref{eq:decompLeb}.
\cqfd

\btheo\label{theo:main2} Assume that $I=J$, that $X$ is a manifold,
that the Gibbs measure $m_{P}$ is finite and mixing under the geodesic
flow on $\Ga\bs \G X$, and that the skinning measure $\sigma^+_{I}$ is
nonzero and finite.  For the narrow convergence of measures on $\Ga\bs
\G X$, we have
\begin{equation}\label{eq:mainequidman}
\lim_{T\ra+\infty}\;\frac{\delta_P\;\|m_P\|}{T\;e^{\frac{\delta_P}{2}T}\;
\|\sigma^+_I\|^2}\;\sum_{[\ga]\in {\mathfrak R}_{II}(T)} 
e^{\frac{1}{2}\int_{[\ga]}P}\Leb_{[\ga]}\;=\frac{1}{\|m_P\|}\;m_P\,.
\end{equation}
Furthermore, if $\wt P=0$, if $X=\wt M$ is a symmetric space, if
$M=\Ga\bs \wt M$ has finite volume and exponentially mixing geodesic
flow, then there exists $k\in\NN$ such that for every compact subset
$K$ of $\Ga \bs\G X$ and every $C^k$-smooth function $\psi:(\Ga \bs\G
X)\ra \CC$ with support in $K$ and $W^{k,2}$-Sobolev norm
$\|\psi\|_{k,2}$, there is an additive error term of the form
$\bigO_K(\frac{\|\psi\|_{k,2}}{T})$ in Equation
\eqref{eq:mainequidman} when evaluated on $\Psi$.
\etheo

As in the case of Theorem \ref{theo:main}, the assumptions for the
error term estimate in Theorem \ref{theo:main2} are satisfied for
instance when $\wt M=\hnr$ and $\Ga$ is geometrically finite, or when
$\Ga$ is an arithmetic lattice.

\medskip
\dem Since $I=J$, the map $\Ga(\alpha,\beta)\mapsto \Ga(\beta,\alpha)$
preserves the multiset $IJ=II$. By Proposition
\ref{prop:bijcomperprecloxo}, this gives, for every $s>0$, an
involution
\[
\rho= \rho_{\alpha,\beta}=\Ga\wt\rho_{\alpha,\beta}\;\mapsto\;
\overline{\rho}= \rho_{\beta,\alpha}=\Ga\wt\rho_{\beta,\alpha}
\]
of $\Perp(\F_I,\F_I,s)$ such that $\mult(\rho)=\mult(\overline{\rho})$
by Equation \eqref{eq:multcomperpsym} and $\int_\rho P=
\int_{\overline{\rho}} P$ by Equation \eqref{eq:reversamplitude}.

Hence, as in order to obtain the left equality in Equation
\eqref{eq:relatNINPerp}, by using Equations \eqref{eq:relatperiod} and
\eqref{eq:decompLeb} (see also Equation \eqref{eq:detailsdecompLeb}),
and since $\|\sigma^-_I\|=\|\sigma^+_I\|$ as mentioned at the end of
Section \ref{sec:dyn}, we have
\begin{align*}
&\frac{\delta_P\;\|m_P\|}{T\;e^{\frac{\delta_P}{2}T}\;
\|\sigma^+_I\|^2}\;\sum_{[\ga]\in {\mathfrak R}_{II}(T)} 
e^{\frac{1}{2}\int_{[\ga]}P}\Leb_{[\ga]}\;\\=\;&
\frac{\delta_P\;\|m_P\|}{T\;e^{\frac{\delta_P}{2}T}\;
\|\sigma^+_I\|^2}\;\sum_{\rho\in \Perp(\F_I,\F_I,\frac{T}{2})} 
e^{\int_{\rho}P}\big(\Leb_{\rho}+\Leb_{\overline{\rho}}\big)
\\=\;&
\frac{\delta_P\;\|m_P\|}{\frac{T}{2}\;e^{\delta_P\frac{T}{2}}\;
\|\sigma^-_I\|\;\|\sigma^+_I\|}\;\sum_{\rho\in \Perp(\F_I,\F_I,\frac{T}{2})} 
e^{\int_{\rho}P}\Leb_{\rho}\,.
\end{align*}
This measure narrow converges to $\frac{1}{\|m_P\|}\;m_P$ as
$\frac{T}{2}\ra+\infty$ by \cite[Theo.~2]{ParPau25b}. Theorem
\ref{theo:main2} follows, with its error term given by the same error
term in \cite[Theo.~1]{ParPau25b}.
\cqfd

\medskip
Using Proposition \ref{prop:bijdihedralrecloxo} and Equation
\eqref{eq:multperiodorbit}, Theorem \ref{theo:main} translates into
counting statements of $\{I,J\}$-dihedral subgroups of $\Ga$ and of
$\{I,J\}$-reversible periodic orbits of the geodesic flow on $\Ga\bs
\G X$, counted with multiplicities and weights. Similarly, Theorem
\ref{theo:main2} translates into an equidistribution statement for the
Lebesgue measures $\Leb_\OOO$ on the $\{I,I\}$-reversible periodic
orbits $\OOO$.  We leave to the willing reader the task to deduce
counting statements for the $\{I,J\}$-primitive $\{I,J\}$-reversible
periodic orbits and for the maximal $\{I,J\}$-dihedral subgroups of
$\Ga$, as well as an equidistribution statement for the Lebesgue
measures $\Leb_\OOO$ on the $\{I,I\}$-primitive $\{I,I\}$-reversible
periodic orbits $\OOO$.

\section{Examples}
\label{ex:plicit}

In this section, we illustrate Theorem \ref{theo:main} with a number
of examples for groups acting on negatively curved symmetric spaces
and on Bruhat-Tits trees. We only consider zero potentials in these
examples, hence we assume $\wt P=0$ throughout Section
\ref{ex:plicit}, and we denote as in the Introduction by
$\N_{I,J}=\N_{I,J,0}$ the counting function (with multiplicities) of
conjugacy classes of $\{I,J\}$-reversible loxodromic elements of $\Ga$.

\subsection{Examples in real hyperbolic spaces}
\label{subsec:realhyp}

Let $\Ga\bs\hnr$ be a finite volume real hyperbolic orbifold of
dimension $n\ge 2$ (with constant sectional curvature $-1$). The
critical exponent of $\Ga$ is $\delta_\Ga=n-1$.  We normalise the
Patterson-Sullivan density $(\mu_{x})_{x\in\hnr}$ such that
$\|\mu_{x}\| = \Vol(\SSS^{n-1})$ for all $x\in\hnr$.   
The total mass of the Bowen-Margulis measure is 
\begin{equation}\label{eq:mbm}
\|m_{\rm BM}\| =2^{n-1}\Vol(\SSS^{n-1})\Vol(\Ga\bs\hnr)\,
\end{equation}
by Assertion (1) of \cite[Prop.~20]{ParPau17ETDS}.  Let $D$ be a
totally geodesic submanifold of $\hnr$ of dimension $k\in\llbracket0,
n-1\rrbracket$. Then $\wt\sigma^+_D=\wt\sigma^-_D= \Vol_{\normalpm D}$
by Assertion (3) of \cite[Prop.~20]{ParPau17ETDS}.\footnote{The case
$k=0$ was not considered in the cited result, but this case is
immediate.}  In particular, with $\Ga_D$ the stabiliser in $\Ga$ of
$D$, if $\Ga_{D}\bs D$ is a properly immersed finite volume
suborbifold of $\Ga\bs\hnr$ and if $\D=(\ga D)_{\ga\in\Ga}$, then
\[
\|\sigma^\pm_\D\|=\Vol(\Ga_{D}\bs\normalpm D)\,.
\]
If $m$ is the number of elements of $\Ga$ that pointwise fix $D$, 
by \cite[Prop.~20 (3)]{ParPau17ETDS}, we have
\begin{equation}\label{eq:skinmass}
\|\sigma^\pm_\D\|=\frac{1}{m}\Vol(\SSS^{n-k-1})\Vol(\Ga_{D}\bs D)\,.
\end{equation}

For every $\alpha\in I_\Ga$, we denote by $\Ga_{F_\alpha}$ the
stabilizer in $\Ga$ of the fixed point set $F_\alpha$ of $\alpha$, and
by $m(\alpha)$ the order of the pointwise stabiliser of $F_\alpha$ in
$\Ga$. Let $I$ and $J$ be two nonempty $\Ga$-invariant subsets of
$I_\Ga$. Let
\begin{equation}\label{eq:defiSigmaI}
\Sigma_I=\sum_{[\alpha]\in \Ga\bs I}\;\frac{1}{m(\alpha)}
\Vol(\SSS^{n-\dim F_\alpha-1})\,\Vol(\Ga_{F_\alpha}\bs F_\alpha)>0\,.
\end{equation}

\bcoro\label{coro:realhyp} Assume that $\Sigma_I$ and $\Sigma_J$
are finite. Then, as $T\ra+\infty$, we have
\[
\N_{I,J}(T)\sim\frac{\Sigma_I\;\Sigma_J}{(n-1)2^{n-1}\Vol(\SSS^{n-1})
  \Vol(\Ga\bs\hnr)}\, e^{\frac{n-1}{2}\;T}\,. 
\]
If $\Ga\bs\hnr$ has finite volume, then there is an additive error
term of the form $\bigO(e^{(\frac{n-1}{2}-\kappa)T})$ for some
$\kappa>0$.
\ecoro

\dem Since $\Ga\bs\hnr$ is locally symmetric with finite volume, the
Bowen-Margulis measure $m_{\rm BM}$ of $\Ga\bs T^1\hnr$ is finite and
mixing. By a summation of Equation \eqref{eq:skinmass} over the
$\Ga$-orbits in $I$ (and recalling that $\sim_I$ is trivial), we have
$\|\sigma^-_I\|=\Sigma_I$. Hence $\|\sigma^-_I\|$ and $\|\sigma^+_J\|$
are nonzero and finite. The asymptotic behavior of $\N_{I,J}(T)$ then
follows from the convergence claim in Theorem \ref{theo:main}
(\hyperlink{theomaincount1}{1}) when $X=\hnr$ and $\wt P=0$ so that
$\delta_P^*=\delta_P=\delta_\Ga=n-1$, and from Equation
\eqref{eq:mbm}.

The above error term follows from the error term in Theorem
\ref{theo:main} (\hyperlink{theomaincount1}{1}), since the geodesic
flow on $\Ga\bs T^1\hnr$ is exponentially mixing by
\cite{LiPan22}. \cqfd

\medskip
We recover \cite[Thm.~1.1]{ErlSou24} in an equivalent form as a
corollary of Corollary \ref{coro:realhyp}, adding 
an error term to their result.

\bcoro\label{coro:sarnak} Let $\Ga$ be a lattice in $\PSL_2(\RR)$ that
contains elements of order $2$. Then there exists $\kappa>0$ such
that, as $T\ra+\infty$, we have
\[
\frac{1}{2}\N_{I_\Ga,I_\Ga}(T)=\frac{\big(\sum_{[\alpha]\in\Ga\bs I_\Ga}
  \frac{1}{|Z_\Ga(\alpha)|}\big)^2}{4\;|\chi^{\rm orb}(\Ga\bs\hdr)|}
\,e^{\frac T2}(1+\bigO(e^{-\kappa T}))\,.
\]
\ecoro

Note that the factor $\frac{1}{2}$ comes from the extra symmetry in
$\N_{I_\Ga,I_\Ga}(T)$ (that might not exists in $\N_{I,J}(T)$ when
$I\cap J=\emptyset$, which is the reason why we chose not to
renormalize by $2$). Otherwise said, $\frac{1}{2}\N_{I_\Ga,I_\Ga}(T)$
counts the number of conjugacy classes of dihedral subgroups $D$ of
$\Ga$, since if $(\alpha,\beta)\in I_\Ga\times I_\Ga$ is such that
$D=\langle \alpha\rangle*\langle\beta\rangle$, then we also have
$D=\langle\beta\rangle*\langle \alpha\rangle$. See for instance Lemma
\ref{lem:computmultprim} which gives the factor $2$ in the primitive
case, considered in \cite[Thm.~1.1]{ErlSou24}. Also note that the
asymptotic growth of the conjugacy classes of dihedral subgroups $D$
and of the conjugacy classes of maximal dihedral subgroups $D$ are the
same.

\medskip
\dem Since $\PSL_2(\RR)$ is the orientation preserving isometry group
of $\hdr$, the elements of $I_\Ga$ act on $\hdr$ by half-turns, each
one fixing a single point. If $\alpha\in I_\Ga$, then the order of the
stabiliser of $F_\alpha$ is $|Z_\Ga(\alpha)|$ (see Lemma
\ref{lem:cardequivmodI} and the comment following its
statement).\footnote{In this case, the centralizer and normalizer of
$\alpha$ are the same groups.}  Note that $\Ga$ contains only finitely
many conjugacy classes of order $2$ elements, since a small enough
Margulis cusp neighborhood in $\Ga\bs\hdr$ contains no nontrivial
orbifold point. Hence $\Ga\bs I_\Ga$ is finite and by Corollary
\ref{coro:realhyp}, for some $\kappa>0$, as $T\ra+\infty$, we have
\[
\frac{1}{2}\N_{I_\Ga,I_\Ga}(T)=\frac{\big(\sum_{[\alpha]\in\Ga\bs I_\Ga}
  \frac{2\pi}{|Z_\Ga(\alpha)|}\big)^2}{4\;(2\pi\Vol(\Ga\bs\hdr))}
\,e^{\frac T2}(1+\bigO(e^{-\kappa T}))\,.
\]
Corollary \ref{coro:sarnak} then follows using the Gauss-Bonnet
formula that relates the hyperbolic area $\Vol(\Ga\bs\hdr) $ of
$\Ga\bs\hdr$ with the orbifold Euler characteristic $\chi^{\rm orb}
(\Ga\bs\hdr)$ of $\Ga\bs\hdr$~: we have $\Vol(\Ga\bs\hdr)
=-2\pi\,\chi^{\rm orb} (\Ga\bs\hdr)$.
\cqfd

\bexer\label{ex:sarnak} Let $p\ge3$ be an integer. Let $\Ga_p$ be the
discrete group of orientation-preserving isometries of $\hdr$
generated by the involution $\alpha: z\mapsto -\frac 1z$, with
$F_\alpha=\{i\}$, and the parabolic element $\ga_p:z\mapsto z+2\cos
\frac\pi p$. It is called the {\em Hecke triangle group} of signature
$(2, p, \infty)$, see for instance \cite[\S 11.3, p.~293]{Beardon83}
or \cite{HaaSer86}. It has a presentation $\langle \alpha, \ga_p\mid
\alpha^2= (\ga_p\alpha)^p=1\rangle$. We have $\Ga_3=\PSL_2(\ZZ)$. The
element $\ga_p\alpha$ is elliptic of order $p$. If $p$ is even, then
$\beta=(\ga_p\alpha)^{\frac p2}$ is an orientation-preserving
involution of $\Ga_p$ with $F_{\beta}=\big\{ e^{i\frac\pi p}\big\}$,
and we have $I_{\Ga_p}=I_\alpha\cup I_\beta$.  The set
\[
\Big\{z\in\CC:-\cos\frac\pi p\le\Re z\le \cos\frac\pi p,\, |z|\ge 1\Big\}
\]
is a fundamental polygon of $\Ga_p$ with boundary identifications
given by $\alpha$ and $\ga_p$. It is easy to check that $\vol(\Ga_p
\bs\hdr)= \pi(1-\frac{2}{p})$, that $\Sigma_{I_\alpha}=\pi$, and when
$p$ is even, that $\Sigma_{I_{\beta}}=\frac{2\pi}p$.  See
\cite{DasGon24} for further details on strongly reversible elements of
Hecke triangle groups. The figure below shows the
$\{I_\alpha,I_\beta\}$-reversible closed geodesics of length at most
$6$ on $\Ga_4\bs\hdr$.

\begin{center}
\includegraphics[width=3.3cm]{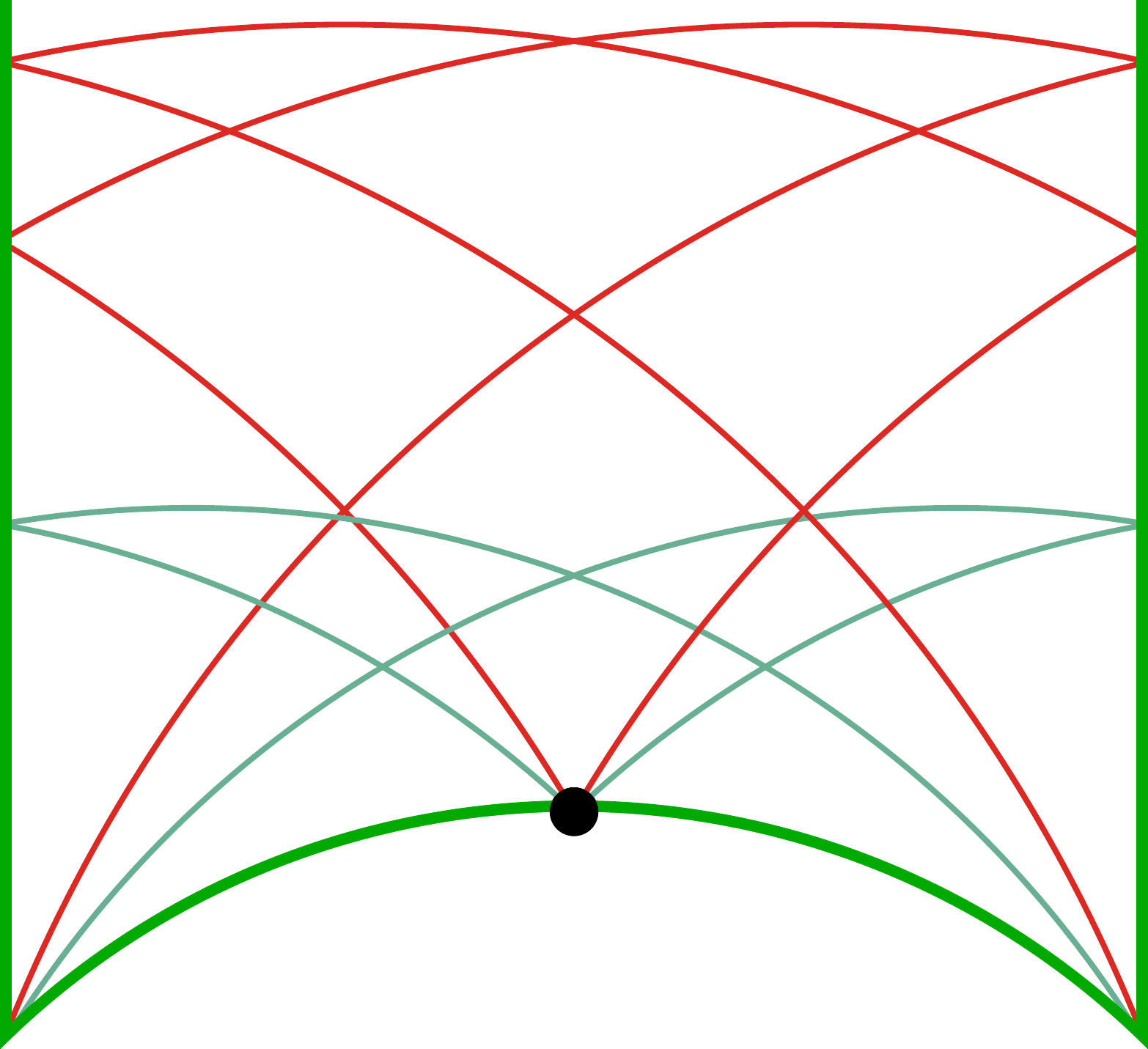}
\end{center}

The figure below shows on the left the
$\{I_\alpha,I_\beta\}$-reversible closed geodesics of length at most
$10$, and on the right the $\{I_\beta,I_\beta\}$-reversible closed
geodesics of length at most $12$ of $\Ga_6\bs\hdr$ in the standard
fundamental polygon of $\Ga_6$.  Note that the
$\{I_\beta,I_\beta\}$-reversible closed geodesic passing through $i$
in the right picture is the square of one of the
$\{I_\alpha,I_\beta\}$-reversible closed geodesics in the left picture
(the blue one), as explained in Remark \ref{rem:defIrecip1}
(\hyperlink{defIrecip1_3}{3}).
\begin{center}
\includegraphics[width=3.3cm]{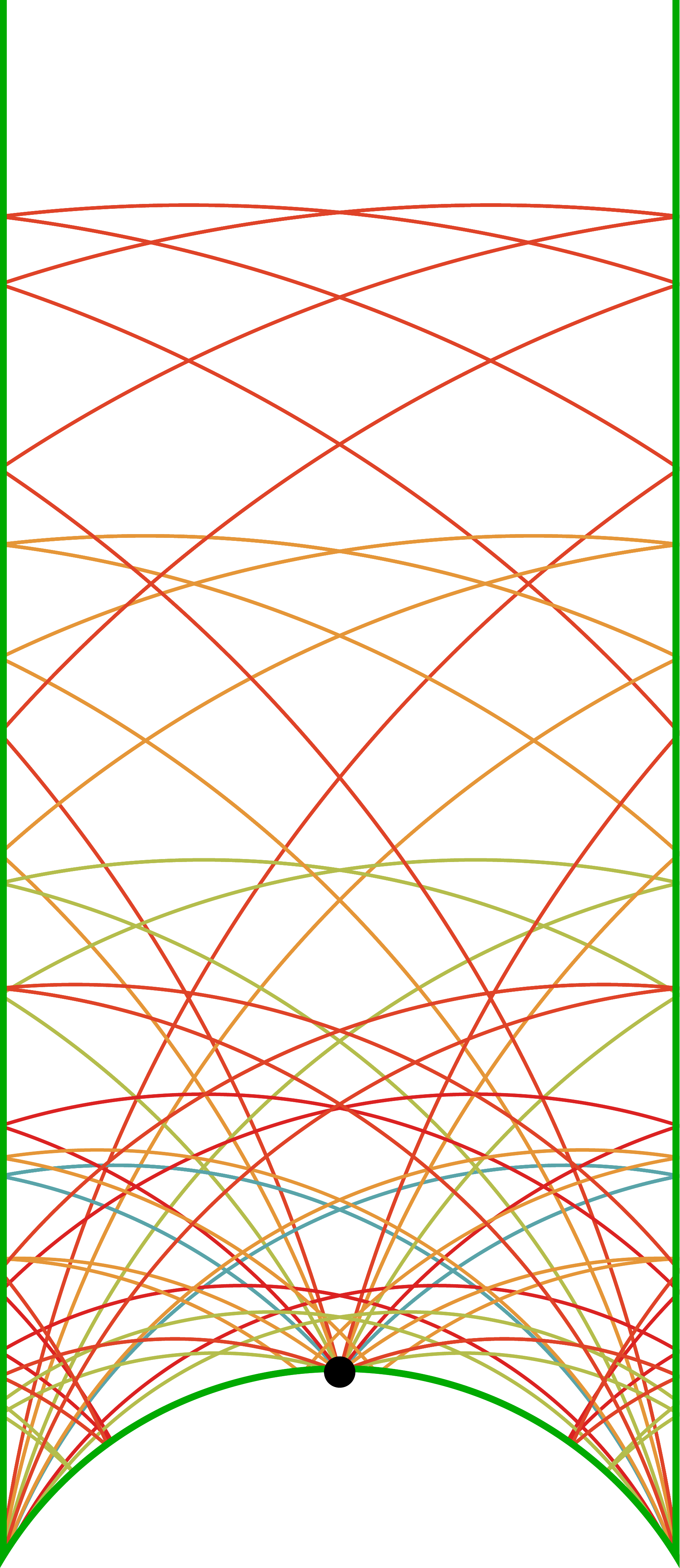}\hspace{2cm}
  \includegraphics[width=3.3cm]{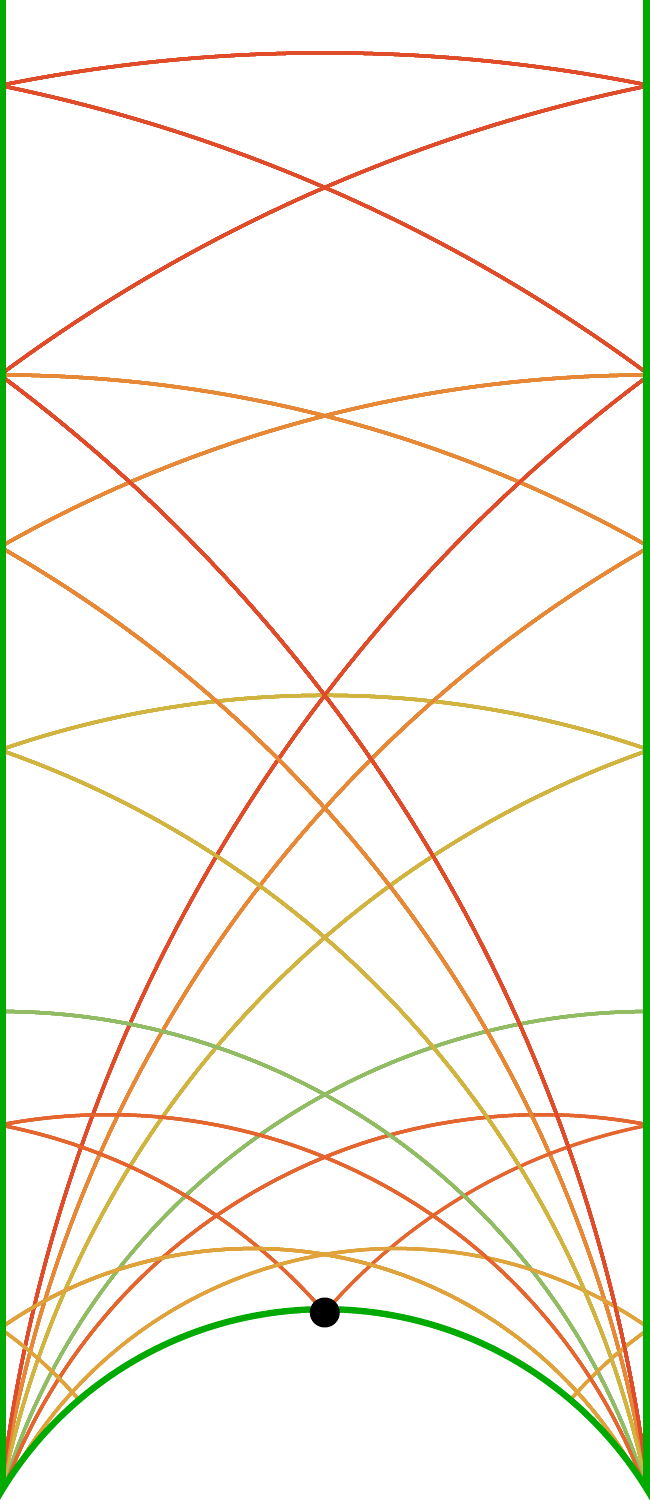}
\end{center}

When $p$ is odd, we have $I_{\Ga_p}=I_{\alpha}$, and Corollary
\ref{coro:realhyp} gives, for some $\kappa>0$,
\[
\N_{I_{\Ga_p},I_{\Ga_p}}(T)=
\frac{1}{4(1-\frac{2}{p})}\;e^{\frac T2}(1+\bigO(e^{-\kappa T})\,.
\]
Since  $\Ga_3=\PSL_2(\ZZ)$,  we get 
\[
\frac{1}{2}\N_{I_{\Ga_3},I_{\Ga_3}}(T)=
\frac{3}{8}\;e^{\frac T2}(1+\bigO(e^{-\kappa T}))\,,
\]
recovering \cite[Thm.~2 (13)]{Sarnak07} albeit with a weaker error
term than loc.~cit.

When $p$ is even, we have $\Sigma_{I_{\Ga_p}}=\Sigma_{I_\alpha}+
\Sigma_{I_\beta} =\pi(1+\frac 2p)$, and Corollary \ref{coro:realhyp}
gives, for some $\kappa>0$,
\[
\N_{I_{\Ga_p},I_{\Ga_p}}(T)=
\frac{(1+\frac 2p)^2}{4(1-\frac{2}{p})}\;e^{\frac T2}(1+\bigO(e^{-\kappa T})\,.
\]
The growth of the number of strongly reversible loxodromic elements of
$\Ga_p$ in terms of word length in the generators is studied in
\cite{DasGon24b}, and in \cite{BasSuz22} for $\Ga_3=\PSL_2(\ZZ)$.
\eexer

Let $P$ be a hyperbolic Coxeter polytope in $\hnr$ with finite nonzero
volume.  Let $S$ be the standard Coxeter generating system consisting
of the reflexions along the codimension $1$ faces of $P$.  Let $W$ be
the nonelementary subgroup of the isometry group of $\hnr$ generated
by $S$, called a {\it hyperbolic Coxeter group}, which is discrete by
Poincaré's theorem. Let $I_S$ be the set of conjugates by the elements
of $W$ of the elements of $S$. Note that in general, we have $I_S\neq
I_W$,\footnote{As mentioned by J.~Parker, for every
hyperbolic Coxeter system $(W,S)$, we have $I_S=I_W$ if and only if
the dihedral angles of $P$ are all odd submultiples of $2\pi$.
} a reason why not to restrict this paper to the case $I=J=I_\Ga$ with
the notation of Sections \ref{sec:proofs} and \ref{sec:eqproof}. The
fixed point set $F_\alpha$ of an element $\alpha\in I_S$ is called a
{\it wall} of $(W,S)$.  All walls are totally geodesic submanifolds of
dimension $n-1$.
  
Finite nonzero volume hyperbolic Coxeter polytopes do not exist if $n>
995$, and are known to exist only in dimensions $m\leq 19$ and
$m=21$. Their classification is known only in dimension $2$ and $3$.
See for instance \cite{FelTum22} for references on hyperbolic Coxeter
systems.

\medskip
\noindent{\bf Proof of Corollary \ref{coro:Coxeter}.}  The Coxeter
polytope $P$ is a fundamental domain for the Coxeter group $W$, and
therefore
\[
\Vol(W\bs\hnr)=\Vol P\,.
\]
Let $\E$ be the set of codimension $1$ faces of $P$. Each $f\in \E$ is
a fundamental domain for the action on the wall $F$ of $(W,S)$
containing $f$ of the stabilizer of $F$ in $W$, so that $\Vol(\Ga_{F}
\bs F)=\Vol(f)$. Furthermore, $\Vol(f)$ is finite unless $n=2$ and $f$
is a geodesic ray or a geodesic line. This has been excluded by the
assumption of Corollary \ref{coro:Coxeter} that if $n=2$, then $P$ is
compact.  The boundary of $P$ is mapped injectively into the quotient
orbifold and the pointwise stabiliser of each wall is generated by the
corresponding reflexion.  The map from $S$ to $W\bs I_S$ which sends
$s\in S$ to $W s$ is bijective. The dimension of the fixed point set
of every element of $I_S$ is $n-1$ and $\Vol(\SSS^0)=2$. Hence with
the definition \eqref{eq:defiSigmaI}, the finiteness assumption of
\[
\Sigma_{I_S}=\sum_{f\in \E}\;\frac{1}{2}\Vol(f)\Vol(\SSS^{n-(n-1)-1})
=\Vol(\partial P)
\]
of Corollary \ref{coro:realhyp} applied with $\Ga=W$ and $I=J=I_S$ is
satisfied. The result then follows from this corollary.
\cqfd

\bexer\label{ex:extendedmodular} Let $p,q,r$ be in $(\NN\ssm\{0,1\})
\cup \{\infty\}$ with $p\leq q\leq r$ and $\frac{1}{p} + \frac{1}{q} +
\frac{1}{r} <1$. Let $W$ be a {\it $(p,q,r)$-hyperbolic triangle
  group}, that is, the hyperbolic Coxeter group with Coxeter polytope a
hyperbolic triangle $P$ in $\hdr$ with angles $\frac{\pi}{p},
\frac{\pi}{q}, \frac{\pi}{r}$. The area of $P$ is
\[
\vol(P)=V_{p,q,r}=\pi-\frac{\pi}{p}-\frac{\pi}{q} -\frac{\pi}{r}\,.
\]
Assume first that $r< \infty$. By the hyperbolic law of cosines, the
perimeter of $P$ is
{\small
\begin{align*} &\vol(\partial P)=L_{p,q,r}=\\&
\arcosh\Big(\frac{\cos\frac{\pi}{r}+\cos\frac{\pi}{p}\cos\frac{\pi}{q}}
{\sin\frac{\pi}{p}\sin\frac{\pi}{q}}\Big)
+\arcosh\Big(\frac{\cos\frac{\pi}{q}+\cos\frac{\pi}{r}\cos\frac{\pi}{p}}
{\sin\frac{\pi}{r}\sin\frac{\pi}{p}}\Big)
+\arcosh\Big(\frac{\cos\frac{\pi}{p}+\cos\frac{\pi}{q}\cos\frac{\pi}{r}}
{\sin\frac{\pi}{q}\sin\frac{\pi}{r}}\Big)\,.
\end{align*}
}

\noindent Hence by Corollary \ref{coro:Coxeter}, with $I_S$ the set of
the conjugates in $W$ of the reflexions along the sides of $P$, as
$T\ra +\infty$, for some $\kappa>0$, we have
\[
\frac{1}{2}\N_{I_S,I_S}(T)=\frac{L_{p,q,r}^2}{8\,\pi\,V_{p,q,r}}
\,e^{\frac{T}{2}}\,(1+\bigO(e^{-\kappa T}))\,.
\]

If $r=\infty$, we have $\vol(\partial P)=+\infty$, and Corollary
\ref{coro:Coxeter} cannot be applied.  By the correspondence in Lemma
\ref{lem:prodinversion} between strongly reversible closed geodesics
and common perpendiculars, replacing, in the proof of Theorem
\ref{theo:main} (\hyperlink{theomaincount1}{1}), the use of
\cite[Coro.~12.3]{BroParPau19} by the use of
\cite[Theo.~6]{ParPau25a}, we can prove that there exists a constant
$c>0$ such that
\[
\N_{I_S,I_S}(T)= c\, T^2\,e^{\frac{T}{2}} +\bigO(Te^{\frac{T}{2}})\,.
\]
\eexer

\bexer\label{ex:Gauss} Let $\operatorname{\text{Möb}}$ be the Möbius
group of the Riemann sphere $\partial\htr=\CC\cup\{\infty\}$, acting
on the upper halfspace model $\CC\times\interval[open]0 {+\infty}$ of
$\htr$ by isometries given by the Poincar\'e extension. Let
$W\subset\operatorname{\text {Möb}}$ be the {\em extended Gaussian
  modular group} generated by the set of reflexions
\[
S=\big\{z\overset{\alpha_1}\mapsto\overline z,\;
z\overset{\alpha_2}\mapsto\overline z+i,\;
z\overset{\alpha_3}\mapsto-\overline z,\;
z\overset{\alpha_4}\mapsto-\overline z+1, \;
z\overset{\alpha_5}\mapsto\frac 1{\overline z}\big\}\,.
\]
The pair $(W,S)$ is a hyperbolic Coxeter system with Coxeter
polyhedron
\[
P=\big\{(z,t)\in\CC\times\interval[open]0{+\infty}:
0\le \Im(z),\Re(z)\le\frac 12\,,\;|z|^2+t^2\ge 1\big\}\,.
\]

Let $\zeta_K$ be Dedekind's zeta function of a quadratic imaginary
number field $K$. Since we have $[\PSL_2(\ZZ[i]):W]=2$, Humbert's
formula proved for example in \cite[Sect.~8.8]{ElsGruMen98} gives
\[
\Vol(W\bs\htr)=\frac{\Vol(\PSL_2(\ZZ[i])\bs\htr)}{[\PSL_2(\ZZ[i]):W]}=
\frac{\frac{2}{\pi^2}\,\zeta_{\QQ(i)}(2)}{[\PSL_2(\ZZ[i]):W]}
=\frac{\zeta_{\QQ(i)}(2)}{\pi^2}\,.
\]

\medskip
\noindent\begin{minipage}{7.3cm} ~~~ In order to use Corollary
\ref{coro:Coxeter}, we compute the area of the boundary of $P$.

~~ The codimension $1$ faces of $P$ contained in the hyperplanes
defined by equations $\Re(z)=0$ and $\Im(z)=0$ are hyperbolic
triangles with angles $0$, $\frac\pi3$ and $\frac\pi2$, each of area
$\frac\pi 6$.  The other two vertical sides are isometric with the
ideal triangle of the real hyperbolic plane $\hdr$ with vertices at
$\infty$, $\frac{\sqrt 3\, i}2$ and $\frac 12+\frac i{\sqrt 2}$ of
area $\frac\pi 2-\arctan(\sqrt 2)$ each. Recall the formula $\cos
c=\cos^2a$ for a right-angled spherical triangle with side lengths
$a,c,a$ and a right angle opposite the side of length $c$. The bottom
of the polyhedron is a Saccheri quadrilateral with three right angles
and the angle $\theta=\arccos(\frac 13)$ at $(\frac{1+i}2, \frac
1{\sqrt 2})\in\htr$, of area $\frac\pi 2-\arccos(\frac 13)$. By
trigonometric angle sum identities, we have $2\arctan(\sqrt 2)+
\arccos(\frac 13)=\pi$.
\end{minipage}
\begin{minipage}{7.6cm}
\begin{center}
\begin{picture}(0,0)%
\includegraphics{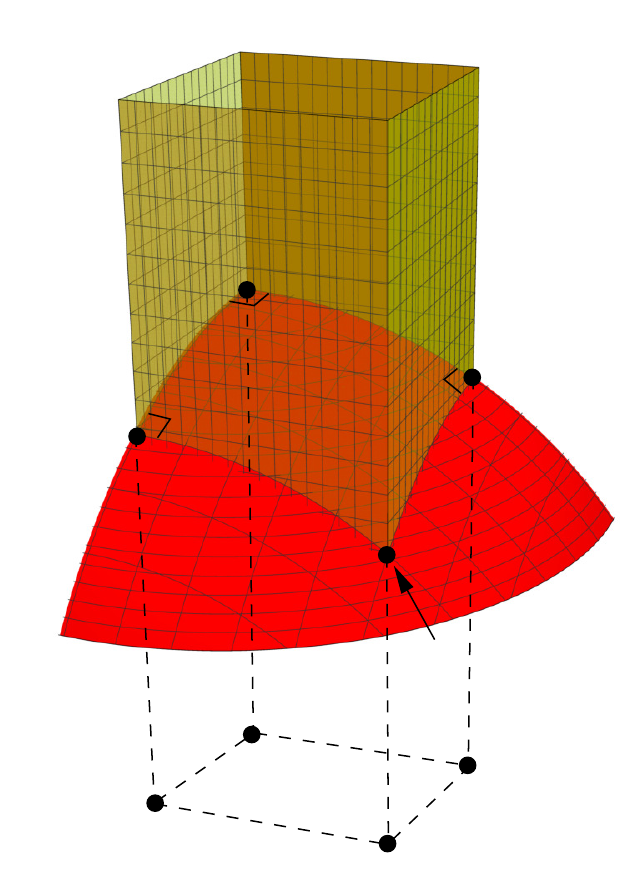}%
\end{picture}%
\setlength{\unitlength}{3812sp}%
\begingroup\makeatletter\ifx\SetFigFont\undefined%
\gdef\SetFigFont#1#2#3#4#5{%
  \reset@font\fontsize{#1}{#2pt}%
  \fontfamily{#3}\fontseries{#4}\fontshape{#5}%
  \selectfont}%
\fi\endgroup%
\begin{picture}(3153,4385)(1938,-4131)
\put(4361,-1581){\makebox(0,0)[lb]{\smash{{\SetFigFont{11}{13.2}{\rmdefault}{\mddefault}{\updefault}{\color[rgb]{0,0,0}$(\frac{i}{2},\frac{\sqrt{3}}{2})$}%
}}}}
\put(1953,-1869){\makebox(0,0)[lb]{\smash{{\SetFigFont{11}{13.2}{\rmdefault}{\mddefault}{\updefault}{\color[rgb]{0,0,0}$(\frac{1}{2},\frac{\sqrt{3}}{2})$}%
}}}}
\put(3866,-3088){\makebox(0,0)[lb]{\smash{{\SetFigFont{11}{13.2}{\rmdefault}{\mddefault}{\updefault}{\color[rgb]{0,0,0}$(\frac{1+i}{2},\frac{1}{\sqrt{2}})$}%
}}}}
\put(3894,-4067){\makebox(0,0)[lb]{\smash{{\SetFigFont{11}{13.2}{\rmdefault}{\mddefault}{\updefault}{\color[rgb]{0,0,0}$\frac{1+i}{2}$}%
}}}}
\put(4327,-3625){\makebox(0,0)[lb]{\smash{{\SetFigFont{11}{13.2}{\rmdefault}{\mddefault}{\updefault}{\color[rgb]{0,0,0}$\frac{i}{2}$}%
}}}}
\put(3237,-3329){\makebox(0,0)[lb]{\smash{{\SetFigFont{11}{13.2}{\rmdefault}{\mddefault}{\updefault}{\color[rgb]{0,0,0}$0$}%
}}}}
\put(2762,-1092){\makebox(0,0)[lb]{\smash{{\SetFigFont{11}{13.2}{\rmdefault}{\mddefault}{\updefault}{\color[rgb]{0,0,0}$(0,1)$}%
}}}}
\put(2522,-3844){\makebox(0,0)[lb]{\smash{{\SetFigFont{11}{13.2}{\rmdefault}{\mddefault}{\updefault}{\color[rgb]{0,0,0}$\frac{1}{2}$}%
}}}}
\end{picture}%

\end{center}
\end{minipage}

\medskip\noindent Thus,
\[
\Vol(\partial P)=2\frac\pi6+2(\frac\pi2-\arctan\sqrt 2)+
\frac\pi2-\arccos\frac 13=\frac56\pi\,.
\]
Corollary \ref{coro:Coxeter} gives, for some $\kappa>0$,
\[
\frac{1}{2}\N_{I_S,I_S}(T)=\frac{(\frac56\pi)^2}{(3-1)\,2^3\,4\pi\,
  \frac{\zeta_{\QQ(i)}(2)}{\pi^2}}\,e^T(1+\bigO(e^{-\kappa T}))=
\frac{25\,\pi^3}{2304\,\zeta_{\QQ(i)}(2)}\,e^T(1+\bigO(e^{-\kappa T}))\,.
\]
\eexer 

\bexer In Example \ref{ex:Gauss}, the collection $I_S$ does not
contain all order $2$ elements of $W$.  For example, the point
symmetry in the point $(0,1)\in\htr$, which is $\alpha= \alpha_5\circ
\alpha_3\circ\alpha_1:z\mapsto-\frac 1{\overline z}$, is not in $I_S$.
Recall that $I_\alpha$ is the set of the conjugates of $\alpha$ by the
elements of $\Ga$. Using the fact that the stabilizer in $W$ of the
point $(0,1)\in\htr$ has order $8$ and by Corollary
\ref{coro:realhyp}, we have, for some $\kappa>0$,
\begin{align*}
\frac{1}{2}\N_{I_\alpha,I_\alpha}(T)&= \frac{(\frac 18\Vol(\SSS^2))^2}
{2^3\;2\;\Vol(\SSS^2)\;\Vol(W\bs\htr)}\,e^T(1+\bigO(e^{-\kappa T}))
\\ &=\frac{\pi^3}{256\;\zeta_{\QQ(i)}(2)}\,e^T(1+\bigO(e^{-\kappa T}))\,.
\end{align*}
\eexer

\subsection{Examples in complex hyperbolic space}
\label{subsec:complexhyp}

Let $\hnc$ be the complex hyperbolic space of (complex) dimension
$n\ge 2$ with sectional curvature $-4\le K\le -1$ (as in
\cite{ParPau17MA}).  Let $\Ga\bs\hnc$ be a finite volume complex
hyperbolic orbifold of dimension $n$.  The critical exponent of $\Ga$
is $\delta_\Ga=2n$ (see for instance \cite[\S 6]{CorIoz99}).
Normalising the Patterson-Sullivan density $(\mu_{x})_{x\in\hnc}$ as
in \cite[\S 4]{ParPau17MA}, the following computations
\eqref{eq:mBMcomphyp} and \eqref{eq:skincomphyp} follow from Lemma 12
of loc.~cit., using its Assertion (iii) for the Bowen-Margulis
measure, and its assertion (vi) for the skinnning measure. We have
\begin{equation}\label{eq:mBMcomphyp}
\|m_{\rm BM}\| =\frac{\pi^n}{2^{2n-3}(n-1)!}\Vol(\Ga\bs\hnc)\,.
\end{equation}
If $D$ is a complex geodesic line in $\hnc$, with stabiliser $\Ga_D$
in $\Ga$, if $\Ga_{D}\bs D$ is a properly immersed finite volume
suborbifold of $\Ga\bs\hnc$, if $\D=(\ga D)_{\ga\in\Ga}$ and if $m$ is
the number of elements of $\Ga$ that pointwise fix $D$, then
\begin{equation}\label{eq:skincomphyp}
\|\sigma^\pm_\D\|=\frac{\pi^{n-1}}{m\,4^{n-2}(n-2)!}\Vol(\Ga_{D}\bs D)\,.
\end{equation}

For every $\alpha\in I_\Ga$, we denote by $\Ga_{F_\alpha}$ the
stabiliser of $F_\alpha$ in $\Ga$ and by $m(\alpha)$ the order of the
pointwise stabiliser of $F_\alpha$ in $\Ga$. Let $I,J$ be nonempty
$\Ga$-invariant subsets of $I_\Ga$. Define
\[
\Sigma'_I=\sum_{[\alpha]\in \Ga\bs I}
\frac{\Vol(\Ga_{F_\alpha}\bs F_\alpha)}{m(\alpha)}\,.
\]
By the above computations and again a summation on the $\Ga$-orbits in
$I$ and $J$, the next result (which has a version valid for every
$n\geq 2$) follows from Theorem \ref{theo:main}
(\hyperlink{theomaincount1}{1}).

\bcoro\label{coro:complexhyp} Assume that $n=2$ and that $\Sigma'_I$
and $\Sigma'_J$ are finite. As $T\ra+\infty$, we have
\[
\N_{I,J}(T)\sim \frac{\Sigma'_I\;\Sigma'_J}
{2\,\Vol(\Ga\bs\hdc)}\, e^{2T}\,.
\]
If $\Ga$ is arithmetic, then there is an additive error term of the
form $\bigO(e^{(2-\kappa)T})$ for some $\kappa>0$. \cqfd
\ecoro

\bexer[Deraux's lattice]\label{ex:Deraux} Deraux's group
$S(2,\sigma_5)$ is the sporadic equilateral triangle group in $\hdc$
generated by three complex reflexions $R_1$, $R_2$ and $R_3$ of order
$2$, that are cyclically permuted by conjugation by an order $3$ complex
reflexion $J$, with parameter $\sigma_5=\operatorname{Tr}(R_1J)
=e^{-\frac{\pi i}{9}} \big(\frac{\sqrt{5}+i\sqrt{3}}{2}\big)$, see
\cite[Section 3.1]{DerParPau21} and
\cite{Deraux06,Deraux24}.\footnote{The group is called $G(4,4,4;5)$ in
\cite{Deraux06}.}

It is an arithmetic lattice by \cite[Theo.~1.2]{DerParPau21}.  The
orbifold Euler characteristic of $S(2,\sigma_5)\bs\hdc$ is $\chi^{\rm
  orb} =\frac{1}{45}$ by \cite[p.~190]{DerParPau21}. By for instance
ibid.~p.~199 or \cite[p.~720]{HerPau96}, this gives $\Vol(S(2
,\sigma_5))\bs\hdc =\frac{1}{2^4}\frac{8\pi^2}{3}\,\chi^{\rm orb}
=\frac{\pi^2}{270}$, taking into account the fact that
\cite{DerParPau21,HerPau96} normalize the sectional curvature to
satisfy $-1\leq K\leq -\frac{1}{4}$ and that the real dimension of
$\hdc$ is $4$.

The complex reflexions $R_1$, $R_2$ and $R_3$ stabilize complex
geodesic lines, called the {\it mirrors} of $S(2,\sigma_5)$, whose
sectional curvature in our normalization is $-4$.  By the last table
page 23 in \cite{Deraux24}, the quotients of the mirrors by their
stabilizers in $S(2,\sigma_5)$ have the same (by the symmetry $J$)
signature $(0;2,2,6,6)$, hence have orbifold Euler characteristic
$-\frac{2}{3}$ . The Gauss-Bonnet formula thus gives an area
$\frac{1}{2^2}(-2\pi)(-\frac{2}{3})=\frac\pi 3$ to each mirror
quotient. The mirrors have pointwise stabilizers of order $2$, again by
the last table on page 23 in \cite{Deraux24}.

Hence with $I$ the set of conjugates of the elements $R_1,R_2,R_3$,
Corollary \ref{coro:complexhyp} says that the number of conjugacy
classes of $\{I,I\}$-reversible loxodromic elements of $S(2,\sigma_5)$
satisfies, for some $\kappa>0$,
\[
\N_{I,I}(T)= \frac{\big(3\frac{1}{2}\frac{\pi}{3}\big)^2}
  {2\,\frac{\pi^2}{270}} e^{2T}(1+\bigO(e^{-\kappa T}))=
\frac{135}{4}\,e^{2T}(1+\bigO(e^{-\kappa T}))\,.
\]
\eexer

\subsection{Examples in trees}
\label{subsec:tree}

Let $\FF_q$ be a finite field of order a positive power $q$ of a prime
number $p$. Let $K$ be a (global) function field over $\FF_q$. Let $v$
be a (normalised discrete) valuation of $K$, let $K_v$ be the
associated completion of $K$, let $\OOO_v=\{x\in K_v:v(x)\geq 0\}$ be
its valuation ring and let $\mmm_v=\{x\in K_v:v(x)> 0\}$ be its
maximal ideal. Let $q_v=q^{\deg v}$ be the order of the residual field
$\OOO_v/\mmm_v$. Let $R_v$ be the affine function ring associated with
$v$. Let $\zeta_K$ be Dedekind's zeta function of $\OOO_K$. For all
these notions and complements, we refer to \cite{Goss96, Rosen02}, as
well as to \cite[\S 14.2]{BroParPau19} whose notation we will follow.

Let $G$ be the locally compact group $\PGL_2(K_v)=\GL_2(K_v)/
(K_v^\times\id)$. We denote by $\begin{bsmallmatrix} a & b \\ c &
  d \end{bsmallmatrix}\in G$ the image in $G$ of $\begin{psmallmatrix}
  a & b \\ c & d \end{psmallmatrix}\in \GL_2(K_v)$.  The group $G$
acts vertex-transitively on the Bruhat-Tits tree $\XX_v$ of
$(\PGL_2,K_v)$, which is a regular tree of degree $q_v+1$, whose
vertices are the homothety classes modulo $K_v^\times$ of the
$\OOO_{v}$-lattices in ${K}_v\times {K}_v$, and whose boundary at
infinity identifies naturally with the projective line $\PP^1(K_v)=
K_v\cup\{\infty\}$. We denote by $*$ the {\em standard basepoint} of
$\XX_v$, which is the homothety class $[\OOO_v\times\OOO_v]$ of the
$\OOO_v$-lattice $\OOO_v\times \OOO_v$ in $K_v\times K_v$, and whose
stabilizer in $G$ is $\PGL_2(\OOO_v)$. See for instance
\cite{Serre83}, as well as \cite[\S 15.1]{BroParPau19} whose notation
we will follow.
 
Let $\Ga=\PGL_2(R_v)$ be the {\it Nagao lattice} in $G$ (see for
instance \cite{Weil70} as well as \cite[\S 15.2]{BroParPau19}) and let
$\alpha=\begin{bsmallmatrix} \ 0&1\\-1& 0 \end{bsmallmatrix}$, which is
an involution of $\Ga$.

\blemm\label{lem:fixsettree} The fixed point set $F_\alpha$ of
$\alpha$ in $\XX_v$ is reduced to $\{*\}$ if $q_v \equiv 3\bmod 4$,
and is the geodesic ray or line in $\XX_v$ with points at infinity the
one or two square roots of $-1$ in $K_v$ otherwise. If $q_v \equiv
3\mod 4$, then $|\Ga_{F_\alpha}/Z_\Ga(\alpha)| >1$.
\elemm

\dem The involution $\alpha$ belongs to $\PGL_2(\OOO_v)$, hence it
fixes the basepoint $*$. For every $n\in\NN\ssm\{0\}$, the sphere
$S(*,n)$ of radius $n$ in $\XX_v$ centered at $*$, which is preserved
by $\alpha$, identifies naturally with $\PP^1(\OOO_v/\mmm_v^n)$, and
in particular with $\PP^1(\FF_{q_v})=\FF_{q_v}\cup\{\infty\}$ if
$n=1$. The involution $\alpha$ acts on $S(*,1)$ by exchanging $\infty$
and $0$, and by sending $z\in \FF_{q_v}\ssm\{0\}$ to $-1/z\in
\FF_{q_v} \ssm\{0\}$. The basepoint $*$ is an isolated fixed point of
$\alpha$ (and hence is the only fixed point of $\alpha$ since
$F_\alpha$ is a subtree) if and only if $\alpha$ has no fixed point on
$S(*,1)$, that is, if and only if the polynomial $X^2+1$ has no root
in the finite field $\FF_{q_v}$.

We claim that this happens if and only if $q_v \equiv 3\mod 4$. If the
characteristic $p$ is equal to $2$, then $X^2+1=(X+1)^2$ has one and
only one root.  Assume now that $p$ is odd. Then $q_v$, which is a
positive power of $p$, is congruent to $1$ or $3$ modulo $4$. Since
$p\neq 2$, a root of $X^2+1$ is a primitive fourth root of unity in
$\FF_{q_v}$. It is well known that a finite field $\FF_{q'}$ of order
$q'$ contains a primitive $n$-th root of unity if and only if $n$
divides $q'-1$, since the multiplicative group $\FF_{q'}^\times$ is a
cyclic group of order $q'-1$. This proves the claim.

If $X^2+1$ has one (when $p=2$) or two (when $p$ is odd and $q_v\equiv
1\mod 4$) roots in $\FF_{q_v}$, then by Hensel's lemma, these one or
two points give the only fixed points of $\alpha$ in
$\PP^1(\OOO_v/\mmm_v^n)$ as well as in $\PP^1(\OOO_v)=\PP^1(K_v)$.

If $q_v \equiv 3\mod 4$, then we have $F_\alpha=\{*\}$ and the
stabilizer in $G$ of $F_\alpha$ is $\PGL_2(\OOO_v)$. Since $R_v\cap
\OOO_v =\FF_q$ by for instance \cite[Eq.~(14.2)]{BroParPau19}, we then
have $\Ga_{F_\alpha}=\PGL_2(\FF_q)$, which has order $q(q^2-1)$.  An
elementary computation, using the fact that $X^2+1$ has no solution in
$\FF_{q}$ if $q \equiv 3\mod 4$, gives that $|Z_\Ga(\alpha)|=2(q+1)$.
Hence $|\Ga_{F_\alpha}/Z_\Ga(\alpha)|= \frac{q(q-1)}{2}$, which is
larger than $1$ since $q\geq 3$. This proves the lemma.
\cqfd

\bcoro\label{cor:Nagao} If $q_v\equiv 3\mod 4$, then as $n\ra+\infty$
with $n\in 4\,\NN$, the number of conjugacy classes of $\{I_\alpha,
I_\alpha\}$-reversible loxodromic elements of $\Ga=\PGL_2(R_v)$ whose
translation length on $\XX_v$ is at most $n$ satisfies, for some
$\kappa>0$,
\[
\N_{I_\alpha,I_\alpha}(n)=
\frac{q_v}{q^2\,(q^2-1)^2\,(q_v-1)\,\zeta_K(-1)}
\;q_v^{\frac n2}(1+\bigO(e^{-\kappa n}))\,. 
\]
\ecoro

\dem We use the basepoint $x_*=*$ in order to define $V_{\rm even}
\XX_v$.  Note that $\Ga$ acts without inversion on $\XX_v$, see
\cite[II.1.3]{Serre83}. Since $\Ga$ is a tree lattice (see for
instance \cite[\S 15.2]{BroParPau19}), the smallest nonempty
$\Ga$-invariant subtree of $\XX_v$ is $\XX_v$ itself, which is regular
of degree $q_v+1\geq 3$.  By \cite[II.1.2.~Corollary]{Serre83} (see
also \cite[page 331]{BroParPau19}), the length spectrum of $\Ga$ is
$L_\Ga=2\ZZ$. We have $\delta_\Ga=\ln q_v$ by Equation (15.8) of
\cite{BroParPau19}.  Let us normalize the Patterson density
$(\mu_x)_{x\in V\XX_v}$ of $\Ga$ to have total mass $\|\mu_x\|=
\frac{q_v+1}{q_v}$ as in Proposition 15.2 (2) ibid.  Then
\[
\|m_{0,\,{\rm even}}\|=\frac{1}{2}\|m_{\rm BM}\|=
\frac{(q_v+1)\zeta_K(-1)}{q_v}
\]
by Proposition 15.3 (1) ibid. The stabilizer of $*$ in $\Ga$ is
$\Ga\cap\PGL_2(\OOO_v)= \PGL_2(\FF_q)$ using Equation (14.2) ibid, and
$*\in V_{\rm even} \XX_v$. Therefore
\[
\|\sigma^\pm_{I_\alpha,\,{\rm even}}\|=\|\sigma^\pm_{I_\alpha}\|=
\frac{\|\mu_*\|}{\card(\PGL_2(\FF_q))}=\frac{q_v+1}{q(q^2-1)q_v}\,.
\]
Theorem \ref{theo:main} (\hyperlink{theomaincount2}{2}) in the case of
trees with $L_\Ga=2\ZZ$, so that $\delta_0^*=1-e^{-2\delta_\Ga}=
\frac{q_v^2-1}{q_v^2}$, gives the claim.
\cqfd

\medskip
If $K=\FF_q(Y)$ is the field of rational fractions over $\FF_q$ with
one indeterminate $Y$ and $v=v_\infty:\frac{P}{Q}\mapsto \deg Q-\deg
P$ for every $P,Q\in \FF _q[Y]$ is the valuation at infinity, then
$K_v=\FF_q((Y^{-1}))$ is the field of formal Laurent series over
$\FF_q$ with indeterminate $Y^{-1}$, $q_v=q$, and $R_v=\FF_q[Y]$. By
for instance \cite[Thm.~5.9]{Rosen02}, we have
\[
\zeta_K(-1)=\frac 1{(q-1)(q^2-1)}\,.
\]
Corollary \ref{cor:introNagao} in the Introduction follows as a
special case of Corollary \ref{cor:Nagao}.

{\small 
\bibliography{../biblio} 
}

\bigskip
{\small
\noindent

\medskip\noindent \begin{tabular}{l} 
Department of Mathematics and Statistics, P.O. Box 35\\ 
40014 University of Jyv\"askyl\"a, FINLAND.\\
{\it e-mail: jouni.t.parkkonen@jyu.fi}
\end{tabular}

\medskip\noindent \begin{tabular}{l}
Laboratoire de mathématique d'Orsay, UMR 8628 CNRS,\\
Universit\'e Paris-Saclay, 91405 ORSAY Cedex, FRANCE\\
{\it e-mail: frederic.paulin@universite-paris-saclay.fr}
\end{tabular}}

\end{document}